\documentclass[a4paper,11 pt]{amsart}
\usepackage{bbm}
\usepackage{mathrsfs}
\usepackage{amsfonts}
\usepackage{latexsym}
\usepackage{arydshln}
\usepackage{lscape}
\usepackage{cases}
\usepackage[arrow,matrix,all]{xy}
\usepackage{amsmath,amssymb,amscd,amsthm,mathrsfs,dsfont,hyperref,amsfonts}

\theoremstyle{plain}\textwidth=31pc \textheight=50pc

\newtheorem{thm}{Theorem}[section]
\newtheorem{cor}[thm]{Corollary}
\newtheorem{lem}[thm]{Lemma}
\newtheorem{prop}[thm]{Proposition}

\newtheorem{rem}[thm]{Remark}

\def\al{\alpha}
\def\bt{\beta}
\def\dt{\delta}
\def\gm{\gamma}

\def\tt{\theta}
\def\vps{\varepsilon}

\def\lmd{\lambda}

\def\bgm{\Gamma}

\def\kk{\mathbbm{k}}
\def\NN{\mathbb{N}}
\def\rr{\mathbbm{r}}

\def\ot{\otimes}
\def\op{\oplus}

\def\se{\leqslant}
\def\le{\geqslant}
\def\lan{\langle}
\def\ran{\rangle}

\def\Hom{\operatorname {Hom}}

\def\Ext{\operatorname {Ext}}
\def\EXT{\operatorname {EXT}}

\def\Im{\operatorname {Im}}
\def\Ker{\operatorname {Ker}}
\def\rad{\operatorname {rad}}

\def\dim{\operatorname {dim}}

\def\H{\operatorname {H}}

\def\cx{\operatorname {cx}}
\def\ad{\operatorname {ad}}
\def\deg{\operatorname {deg}}

\def\mc{\mathcal}
\def\mk{\mathfrak}
\def\mb{\mathbb}
\def\mbb{\mathbbm}

\def\it{\textit}
\def\t{\text}

\def\ra{\rightarrow}
\def\xra{\xrightarrow}

\title{Ext algebra of  Nichols algebras of type $A_2$}

\author{Xiaolan YU}
\address {Xiaolan YU%\newline \indent Department of Mathematics, Hangzhou Normal University, Hangzhou, Zhejiang 310036, China
\newline Department WNI, University of Hasselt, Universitaire Campus, 3590
Diepeenbeek,Belgium } \email{xiaolan.yu@uhasselt.be }

\author{Yinhuo ZHANG}
\address {Yinhuo ZHANG\newline Department WNI, University of Hasselt, Universitaire Campus, 3590 Diepeenbeek,Belgium } \email{yinhuo.zhang@uhasselt.be}

\begin{document}

\begin{abstract}
We give the full structure of the Ext algebra of a Nichols algebra
of type $A_2$ by using the Hochschild-Serre spectral sequence. As an
application, we show that the pointed Hopf algebras $u(\mathcal{D},
\lmd, \mu)$ with Dynkin diagrams of type $A$, $D$, or $E$, except
for $A_1$ and $A_1\times A_1$ with the order $N_{J}>2$ for at least
one component $J$, are wild.
\end{abstract}

%%% ----------------------------------------------------------------------
\maketitle
%%% ----------------------------------------------------------------------

\section*{Introduction}
For an algebra $R$ over a field $\kk$, its homological properties,
such as the Calabi-Yau property \cite{g2}, AS-regularity
\cite{lpwz}, support varieties \cite{ss}, etc. rely exclusively on
the structure of its Ext algebra $\Ext_{R}^*(\kk,\kk)$.

Nichols algebras play an important role in the classification of
pointed Hopf algebras \cite{ahs,as,as2,fg}. They are braided Hopf
algebras in certain braided monoidal categories.  In \cite{as2}, the
authors showed that if $H$ is a finite dimensional pointed Hopf
algebra  such that its coradical is an abelian group with order not
divisible by primes less than 11, then $H$ is isomorphic to a
deformation of the bosonization of a  Nichols algebra of finite
Cartan type. Thus the study of Nichols algebras not only helps us to
classify pointed Hopf algebras, but also helps us to understand more
about the properties of pointed Hopf algebras. In two recent papers
\cite{fw,pw} support varieties of modules over Hopf algebras are
introduced. It turns out that support varieties are useful tools to
study homological properties and representations of (braided) Hopf
algebras. To define and to compute support varieties over a
(braided) Hopf algebra we need first to understand the Ext algebra
of the (braided) Hopf algebra.  In \cite{a}, the author raised the
question of when the Ext algebra  of a Nichols algebra is still a
Nichols algebra. These facts motivate us to study the structure of
the Ext algebra of a Nichols algebra in this paper.

As a first attempt to explore the structure of  the Ext algebras for
further study,  we will give the full structure of the Ext algebra
of a Nichols algebra of type $A_2$.  First we use the
Hochschild-Serre spectral sequence to get a basis of the  Ext
algebra. We then construct the first segment of the minimal
projective resolution of $\kk$ and give the relations that hold in
the Ext algebra. We calculate the dimensions to verify that these
relations are complete (see Theorems \ref{1} and \ref{2}) . The
relations are braided commutative, which coincides with what have
been proved in \cite{mpsw}, where the authors also showed that the
cohomology ring of a finite dimensional pointed Hopf algebra of
finite Cartan type is finitely generated. Having the generators and
relations of the Ext algebra, we can show that the Ext algebra of a
Nichols algebra is not a Nichols algebra in general (see Proposition
2.15). However, the quotient algebra of the Ext algebra modulo the
ideal generated the nilpotents can be a Nichols algebra (see
Proposition 2.16). This partially answers one of the questions
raised in \cite[Sec. 2.1]{a}.

Finite dimensional pointed Hopf algebras with  abelian group
coradicals have support varieties \cite{fw,mpsw}.  For a pointed
Hopf algebra $A$ of type $A_2$, the support variety of $\kk$ over
$A$ is isomorphic to the variety of $\kk$ over the associated graded
algebra with respect to a certain filtration of $A$. This can be
showed by using the full structure of the Ext algebra of the Nichols
algebra of type $A_2$. Finally, we apply our main results to show
that in many cases, the pointed Hopf algebras $u(\mathcal{D}, \lmd,
\mu)$ constructed in \cite{as2} are wild (Proposition \ref{repre}).

\subsection*{Acknowledgement} The first author deeply thanks Dr.
Jiwei He for helpful discussions and useful comments.

\section{Preliminaries and Notations}

Throughout the  paper, we fix an algebraically closed field $\kk$
with $\t{char}\kk\neq2$. All algebras are assumed to be finite
dimensional and all modules are assumed to be finitely generated
unless otherwise stated.

\subsection{Nichols Algebras and pointed Hopf algebra of Cartan
type} In \cite{as2}, the authors classified finite dimensional
pointed Hopf algebras whose coradicals are  abelian groups.  We need
the following terminology:
\begin{itemize}
\item  a finite abelian group $\bgm$; %In this paper, we always assume that  $\bgm$ is finitely generated;
\item  a Cartan matrix $(a_{ij})\in \mathbb{Z}^{\tt\times \tt}$ of finite type, where $\tt\in \NN$;
\item a set $\mathcal {X}$ of connected components of the Dynkin diagram corresponding
to the Cartan matrix $(a_{ij})$. If $1\se i, j\se \tt$, then $i\sim
j$ means that they belong to the same connected component;
\item elements $g_1,\cdots , g_\tt\in \bgm$ and characters $\chi_1,\cdots, \chi_\tt\in \widehat{\bgm}$ such that
\begin{equation}\label{q}\chi_j(g_i)\chi_i(g_j)=\chi_i(g_i)^{a_{ij}}, \t{  } \chi_i(g_i)\neq 1, \t{   for all $1\se i,j\se \tt$}.\end{equation}
\end{itemize}

The collection $\mc{D}(\bgm, (g_i)_{1\se i\se \tt},(\chi_i)_{1\se
i\se \tt},(a_{ij})_{1\se i,j\se \tt} )$ is called a \textit{datum of
finite Cartan type} for $\bgm$.

For simplicity, we define $q_{ij}=\chi_j(g_i)$. Then equation
(\ref{q}) reads as
\begin{equation}\label{q1}q_{ij}q_{ji}=q_{ii}^{a_{ij}},\;\; q_{ii}\neq 1, \t{   for all $1\se i,j\se \tt$}.\end{equation}

From now on, we assumt that for
 all $1\se i\se \tt$,
\begin{equation}\label{odd}\begin{array}{l}q_{ii} \t{ has odd order, and}\\
\t{the order of } q_{ii} \t{ is prime to }3, \t{ if $i$ lies in a
component $G_2$}.
\end{array}\end{equation}
Since $q_{ij}q_{ji}=q_{ii}^{a_{ij}}$, $1\se i,j\se \tt$, the order
of $q_{ii}$ is constant in each component $J\in \mc{X}$ of the
Dynkin diagram. Let $N_J$ denote this common order.

Given a datum $\mc{D}$, we define a braided vector space as follows.
Let $V$ be a Yetter-Drinfeld module over the group algebra $\kk
\bgm$ with basis  $x_i\in V^{\chi_i}_{g_i}$, $1\se i\se \tt$. Then
$V$ is a braided vector space of diagonal type whose braiding is
given by \begin{equation}\label{bvs}c(x_i\ot x_j)=q_{ij}x_j\ot
x_i,\;\;1\se i,j\se \tt.\end{equation}

Let $\lambda=(\lambda_{ij})_{1\se i<j\se n,i\nsim j}$ be a set of
scalars, such that
$$\lambda_{ij}=0\;\;  \text{if $g_ig_j=1$ or $\chi_i\chi_j\neq
\varepsilon$},$$ where $\varepsilon$ is the identity in
$\widehat{\bgm}$. The set of scalars $\lmd=(\lambda_{ij})_{1\se
i,j\se n,i\nsim j}$ are  called \it{linking parameters}. The algebra
$U(\mc{D},\lmd)$ is defined to be the  quotient Hopf algebra of the
smash product $\kk\lan x_1,\cdots, x_\tt\ran\#\kk \bgm$ modulo the
ideal generated by the following relations
{\small$$\begin{array}{lll}
\t{(Serre relations)}&(\ad_cx_i)^{1-a_{ij}}(x_j)=0,&1\se i,j\se \tt,\;\;i\neq j,\;\;i\sim j,\\
\t{(linking
relations)}&x_ix_j-\chi_j(g_i)x_jx_i=\lmd_{ij}(1-g_ig_j),&1\se
i<j\se \tt,\;\;i\nsim j,
\end{array}$$}
where $\ad_c$ is the braided adjoint representation defined in
\cite[Sec. 1.4]{as}.
%If $x$ is a primitive element, then the braided
%adjoint representation of $x$ is just
%$$\ad_c(x)(y)=m(\id-c)(x\ot y):=[x,y]_c.$$
%$[x,y]_c$ is called a \textit{braided commutator}.

Let $\Phi$ be the root system corresponding to the Cartan matrix
$(a_{ij})$ with $\Pi=\{\al_1,\cdots, \al_\tt\}$ a set of fixed
simple roots. Let $\Phi_J$, $J\in \mc{X}$, be the root system of the
component $J$. Assume that $\mathcal {W}$ is the Weyl group of the
root system $\Phi$. We fix a reduced decomposition of the longest
element
$$w_0=s_{i_1}\cdots s_{i_p}$$ of $  \mc{W}$ as a product of simple
reflections. Then the positive roots $\Phi^+$ are precisely the
followings
$$\bt_1=\al_{i_1}, \;\;\bt_2=s_{i_1}(\al_{i_2}),\cdots, \bt_p=s_{i_1}\cdots s_{i_{p-1}}(\al_{i_p}).$$ If
$\bt_i=\sum_{j=1}^{\tt} m_j\al_j$, then we define
$$g_{_{\bt_i}}=g_1^{m_1}\cdots g_\tt^{m_\tt} \t{ and }\chi_{_{\bt_i}}={\chi}_1^{m_1}\cdots {\chi}_\tt^{m_\tt}.$$ Similarly,
we write $q_{_{\bt_j\bt_i}}=\chi_{_{\bt_i}}(g_{_{\bt_j}})$.

Let $x_{\bt_j}$, $1\se j \se p$, be the root vectors as defined in
\cite[Sce. 2.1]{as2}. Let $(\mu_{\al})_{\al\in \Phi^+}$ be a set of
scalars, such that
$$\mu_{\al}=0 \;\; \text{if $g_{\al}^{N_J}=1$ or $\chi_{\al}^{N_J}\neq \varepsilon$}, \;\;\al\in\Phi^+_J,J\in
\mc{X}.$$ This set of scalars are called \it{root vector
parameters}. The finite dimensional Hopf algebra
$u(\mathcal{D},\lmd, \mu)$ is the quotient of $U(\mathcal{D},\lmd)$
modulo the ideal generated by

{\small$\begin{array}{lll} \t{(root vector
relations)}&\hspace{1cm}x_{\al}^{N_J}-u_{\al}(\mu),&\hspace{1cm}\al\in\Phi^+_J,J\in
\mc{X},
\end{array}$}

where $u_{\al}(\mu)\in \kk \bgm$ is defined inductively on $\Phi^+$
as in \cite[Sec 4.2]{as2}.

Let $V$ be the braided vector space defined as in (\ref{bvs}). The
Nichols algebra $\mathcal {B}(V)$ associated to $V$ is a braided
Hopf algebra in the category of Yetter-Drinfeld modules over
$\kk\bgm$. By \cite[Thm. 5.1]{as2}, it is generated by
$x_1,\cdots,x_\tt$ subject to relations
$$(\ad_cx_i)^{1-a_{ij}}(x_j)=0, \;\;1\se i,j\se \tt,\;\;i\neq j,$$
$$x_{\al}^{N_J}=0, \;\;\al\in\Phi^+_J,J\in \mc{X}.$$
The details about Nichols algebras can be found in \cite{as}.

Corollary 5.2 in \cite{as2} showed that the associated graded Hopf
algebra $\t{Gr} u(\mathcal{D},\lmd,\mu)$ of the algebra
$u(\mathcal{D},\lmd,\mu)$ with respect to the coradical filtration
is $u(\mathcal{D},0,0)$. Moreover, we have that $\mathcal {B}(V)\#
\kk \bgm\cong u(\mathcal{D},0, 0)$. More detailed discussion about
the algebras $U(\mathcal{D},\lmd)$ and $u(\mathcal{D},\lmd,\mu)$ can
be found in \cite{as2}.

The following set
$$\{x_{_{\bt_1}}^{a_1}\cdots x_{_{\bt_p}}^{a_p}\mid 1\se a_i<
N_J,\;\;\;\bt_i\in \Phi^+_J,\;\;\;1\se i\se p\}$$ forms a PBW basis
of the Nichols algebra $\mathcal {B}(V)$ \cite{as2}. As in
\cite[Sec. 2]{mpsw}, define a degree on each element as
$$\deg x_{_{\bt_1}}^{a_1}\cdots x_{_{\bt_p}}^{a_p}=(\sum a_i ht(\bt_i),a_p,\cdots, a_1) \in \NN^{p+1},$$ where $ht(\bt_i)$ is the
height of the positive root $\bt_i$.  That is, if
$\bt_i=\sum_{j=1}^\tt m_j\al_j$, then $ht(\bt_i)=\sum_{j=1}^\tt
m_j$. Order the elements in $\NN^{p+1}$ as follows
\begin{equation}\label{ordering}\begin{array}{l}
(a_{p+1},a_p, \cdots,a_1)<(b_{p+1},b_p, \cdots,b_1)\t{ if and only
if there is some } \\1\se k\se p+1 ,\t{ such that } a_i=b_i \t{ for
}i\le k \t{ and }a_{k+1}<b_{k+1}.\end{array}
\end{equation}

By \cite[Thm. 9.3]{dp}, similar to the proof of Lemma 2.4 in
\cite{mpsw}, we obtain the following lemma.

\begin{lem}\label{nichols filt}
In the Nichols algebra $\mc{B}(V)$, for $j>i$, we have
\begin{equation}\label{comm1}[x_{_{\bt_i}}, x_{_{\bt_j}}]_c=\sum_{\textbf{a}\in \mathbb{N}^p}\rho_{\textbf{a}}x_{_{\bt_1}}^{a_1}\cdots x_{_{\bt_p}}^{a_p}, \end{equation}
where $\rho_{\textbf{a}}\in \kk$ and $\rho_{\textbf{a}}\neq0$ only
when $\textbf{a}=(a_1,\cdots,a_p)$ satisfies that $a_k=0$ for $k\se
i$ and $k\le j$.
\end{lem}

Therefore, if we order PBW basis elements by degree as in
(\ref{ordering}), we obtain a filtration on the Nichols algebra
$\mathcal {B}(V)$. The associated graded algebra
$\mb{G}\mbb{r}\mathcal {B}(V)$ is generated by the root vectors
$x_{_{\bt_{i}}}$, $1\se i\se p$, subject to the relations
$$[x_{_{\bt_{i}}},x_{_{\bt_{j}}}]_c=0, \;\;\text{for all $i<j$};$$
$$x_{_{\bt_{i}}}^{N_J}=0,\;\;\bt_i\in \Phi_J^+,\;\;1\se i\se p.$$

\subsection{Complexity and varieties} We follow the definitions and the notations in \cite{fw}. Let $A$ be a finite
dimensional Hopf algebra and $\H^*(A,\kk):=\Ext_A^*(\kk,\kk)$. The
vector space $\H^*(A,\kk)$ is an associative graded algebra under
the Yoneda product. The subalgebra $\H^{ev}(A,\kk)$ of $\H^*(A,\kk)$
is defined as
$$\H^{ev}(A,\kk)=\op_{n=0}^\infty\H^{2n}(A,\kk).$$ The algebra $\H^{ev}(A,\kk)$ is commutative, since
$\H^*(A,\kk)$ is graded commutative. In the following, we say that a
Hopf algebra $A$ satisfies the \textrm{\textit{assumption
}(\textbf{fg})} if the following conditions hold:
$$\begin{array}{ll}
\t{{\bf(fg1)}}&\t{The algebra $\H^{ev}(A,\kk)$ is finitely
generated.}\\
\t{{\bf(fg2)}}&\t{The $\H^{ev}(A,\kk)$-module $\Ext^*_A(M,N)$ is
finitely generated for any two} \\
&\t{finite dimensional $A$-modules $M$ and $N$.}
\end{array}$$
Under the assumption ({\bf fg}),  the \textit{variety}
$\mathcal{V}_A(M,N)$ for $A$-modules $M$ and $N$ is defined as
$$\mathcal{V}_A(M,N):=\text{MaxSpec}({\H^{ev}(A,\kk)}/I(M,N)),$$ where $I(M,N)$ is the annihilator of the action of $\H^{ev}(A,\kk)$ on $\Ext_A^*(M,N)$. It is a homogeneous ideal of $\H^{ev}(A,\kk)$.  The \textit{support variety} of $M$ is defined as $\mathcal{V}_A(M)=\mathcal{V}_A(M,M)$. By \cite[Thm 6.3]{mpsw}, a
finite dimensional pointed Hopf algebra of the form $u(\mc{D},\lmd,
\mu)$ satisfies the assumption ({\bf fg}).

For a graded vector space $V^\bullet=\op_{n\in\mathbb{Z}^{\le 0}}
V^n$, the growth rate $\gamma(V^\bullet)$ is defined as
$$\gamma(V^\bullet)=min\{c\in \mathbb{Z}, c\le 0\mid \exists b\in
\mathbb{R}, \text{such that} \;\dim V^n\se bn^{c-1}, \text{for all
$n\le 0$}\}.$$

Let $M$ be an $A$-module and $P_*:\;\cdots \ra P_1\ra P_0\ra M\ra 0$
 a minimal projective resolution of $M$. Then the growth rate
$\gamma(P_*)$ is defined to be the \it{complexity} $\cx_A(M)$ of
$M$.

\section{Ext algebras }

It is clear that each Nichols algebra can be written as a twisted
tensor product of a set of Nichols algebras, such that each of them
satisfies that the Dynkin diagram associated to the Cartan matrix is
connected. In \cite{bo}, the authors showed that the $\Ext$ algebra
of a twisted tensor algebra is essentially the twisted tensor
algebra of the $\Ext$ algebras. Therefore, we only need to discuss
the case where the Dynkin diagram is connected. Now we calculate
 the $\Ext$ algebra of a Nichols algebra of type $A_2$.

Let $N$ be an integer, and let $\bar{q}$ be a primitive root of 1 of
order $N$. Let $q_{ij}$, $1\se i,j\se 2$ be roots of 1, such that
$$q_{_{11}}=q_{_{22}}=\bar{q},\;\;\;q_{_{12}}q_{_{21}}=\bar{q}^{-1}.$$

Let $V$ be a 2-dimensional vector space with basis $x_1$ and $x_2$,
whose braiding is given by
$$c(x_i\ot x_j)=q_{ij}x_j\ot x_i,\;\;1\se i,j\se 2.$$
Then $V$ is a braided vector space of type $A_2$.

\subsection{Case $N=2$} As discussed in \cite{ad}, the Nichols algebra
$R=\mathcal {B}(V)$ is isomorphic to the algebra  generated by $x_1$
and $x_2$, with relations
$$x_1x_2x_1x_2+x_2x_1x_2x_1=0,\;\;x_1^2=x_2^2=0.$$ The dimension of
$R$ is 8.

Its $\Ext$ algebra can be calculated directly via the minimal
projective resolution of $\kk$.

Throughout, for an algebra $R$, we write elements in the free module
$R^n$, $n\le 1$, as row vectors. A morphism $f:R^m\ra R^n$ is
described by an $m\times n$ matrix.

\begin{prop}\label{N=2}
Let $R=\mc{B}(V)$ be the algebra mentioned before, then the algebra
$\Ext_R^*(\kk,\kk)$ is generated by $\mk{a}_1$, $\mk{a}_2$ and
$\mk{b}$ with $\deg \mk{a}_1=\deg \mk{a}_2=1$ and $\deg \mk{b}=2$,
subject to the relations
$$\mk{a}_2\mk{a}_1=\mk{a}_1\mk{a}_2=0,\;\;\mk{a}_1\mk{b}=\mk{b}\mk{a}_1,\;\;\mk{a}_2\mk{b}=\mk{b}\mk{a}_2.$$
\end{prop}
\proof We claim that the following complex  is the minimal
projective resolution of $\kk$. \begin{equation}\label{c}\cdots
P_n\xra{d_n} P_{n-1}\ra\cdots P_2\xra{d_2} P_1\xra{d_1} P_0\ra
\kk,\end{equation} where $P_n=R^{n+1}$ and $d_n$ is defined as

$$d_n={\tiny\left(\begin{array} {cccccccc}x_1\\
x_2x_1x_2&x_1\\
&\cdots\\
&&x_2x_1x_2&x_1\\
&&&x_2&x_1x_2x_1\\
&&&&\cdots \\
&&&&&&x_2&x_1x_2x_1\\
&&&&&&&x_2
\end{array}\right)},$$when $n$ is odd and

$$d_n={\tiny\left(\begin{array} {cccccccc}x_1\\
x_2x_1x_2&x_1\\
&\cdots \\
&&x_2x_1x_2&x_1\\
&&&x_2x_1x_2&x_1x_2x_1\\
&&&&x_2&x_1x_2x_1\\
&&&&&\cdots \\
&&&&&&x_2&x_1x_2x_1\\
&&&&&&&x_2
\end{array}\right)}, $$when $n$ is even. Especially, $d_1={\tiny\left(\begin{array} {c}x_1\\
x_2
\end{array}\right)}$.
It is routine to check that (\ref{c}) is indeed a complex.  Now we
use induction to prove the exactness. It is clear that the minimal
projective resolution starts as
$$R^3\xra{d_2}R^2\xra{d_1} R\ra \kk\ra 0,$$ where $d_1={\small\left(\begin{array} {c}x_1\\
x_2
\end{array}\right)}$  and $d_2={\small\left(\begin{array} {cc}x_1&\\
x_2x_1x_2&x_1x_2x_1\\
&x_2
\end{array}\right)}$. Assume that the complex (\ref{c}) is exact
up to  $P_n$. If $n$ is odd, then {\small
$$\begin{array}{cl}&\dim(\Ker
d_n)\\=&(1+\dim P_1+\dim P_3+\cdots\dim P_n)-(\dim P_0+\dim P_2+\cdots\dim P_{n-1})\\
=&4n+5.\end{array}$$} Since the dimension of $R$ is small, we can
calculate the dimension of the submodule $\Im d_{n+1}$ of $P_n$
directly, it is also $4n+5$. Then the complex is exact at $P_{n+1}$.
If $n$ is even, by a similar discussion, we can also conclude that
the complex is exact at $P_{n+1}$.  In this case $\dim(\Ker
d_n)=4n+7$. We have that  $\Im d_n\subseteq \rad P_{n-1}$ for each
$i\le 0$. Therefore, the complex (\ref{c}) is the minimal projective
resolution of $\kk$. Since $\kk$ is a simple module, we have
\begin{equation}\label{equ
mini}\Hom_R(P_n,\kk)\cong\Ext^n_R(\kk,\kk)\end{equation} as vector
spaces for each $n\le 0$. Let
$\mk{a}_1$,$\mk{a}_2\in\Hom_R(P_1,\kk)$ be the functions dual to
$(1,0)$ and $(0,1)$ respectively and $\mk{b}\in\Hom_R(P_2,\kk)$ be
the function dual to $(0,1,0)$.

Let $f_i$, $g_i$ and $h_i$ be the morphisms described by the
following matrices:

$$f_1={\tiny\left(\begin{array}{c}1\\0\end{array}\right)},\;\;\; f_2={\tiny\left(\begin{array} {cc}1&0\\
0&x_2x_1\\
0&0
\end{array}\right)}, \;\;\;f_3={\tiny\left(\begin{array} {ccc}1&0&0\\
0&1&0\\
0&0&0\\
0&0&0\\
\end{array}\right)},$$

$$g_1={\tiny\left(\begin{array}{c}0\\1\end{array}\right)}, \;\;\; g_2={\tiny\left(\begin{array} {cc}0&0\\
x_1x_2&0\\
0&1
\end{array}\right)},\;\;\; g_3={\tiny\left(\begin{array} {ccc}0&0&0\\
0&0&0\\
0&1&0\\
0&0&1\\
\end{array}\right)},$$

$$h_2={\tiny\left(\begin{array} {c}0\\
1\\
0
\end{array}\right)}, \;\;\;h_3={\tiny\left(\begin{array} {cc}0&0\\
1&0\\
0&1\\
0&0\\
\end{array}\right)}.$$
Then we have the following commutative diagrams: {$$\xymatrix{
  P_3\ar[r]^{d_3}\ar[d]_{f_3}&P_2\ar[r]^{d_2}\ar[d]_{f_2}&P_1\ar[r]^{d_1}\ar[d]_{f_1}\ar[dr]&P_0\ar[r]&\kk\\
   P_2\ar[r]^{d_2}&P_1\ar[r]^{d_1}&P_0\ar[r]&\kk},$$}
{$$\xymatrix{
  P_3\ar[r]^{d_3}\ar[d]_{g_3}&P_2\ar[r]^{d_2}\ar[d]_{g_2}&P_1\ar[r]^{d_1}\ar[d]_{g_1}\ar[dr]&P_0\ar[r]&\kk\\
   P_2\ar[r]^{d_2}&P_1\ar[r]^{d_1}&P_0\ar[r]&\kk},$$}
   {$$\xymatrix{
  P_3\ar[r]^{d_3}\ar[d]_{h_3}&P_2\ar[r]^{d_2}\ar[d]_{h_2}\ar[dr]&P_1\ar[r]^{d_1}&P_0\ar[r]&\kk\\
   P_1\ar[r]^{d_2}&P_0\ar[r]^{d_1}&\kk}.$$}
These commutative diagrams show that the relation listed in the
proposition hold.

Let $U$ be the algebra  generated by $\mk{a}_1$, $\mk{a}_2$ and
$\mk{b}$  subject to the relations listed in the proposition. When
$n$ is odd, $U_n$ has a basis
$$\{\mk{a}_1^n, \mk{a}_1^{n-2}\mk{b},\cdots,\mk{a}_1\mk{b}^{\frac{n-1}{2}}, \mk{a}_2\mk{b}^{\frac{n-1}{2}},\cdots,\mk{a}_2^{n-2}\mk{b}, \mk{a}_2^n\}$$ and when $n$ is even, $U_n$ has a basis
$$\{\mk{a}_1^n, \mk{a}_1^{n-2}\mk{b},\cdots,\mk{a}_1\mk{b}^{\frac{n}{2}-1},\mk{b}^{\frac{n}{2}}, \mk{a}_2\mk{b}^{\frac{n}{2}-1},\cdots, \mk{a}_2^{n-2}\mk{b},\mk{a}_2^n\}.$$ They are functions dual to $(1,0\cdots,0)$, $\cdots$, $(0,\cdots,0,1)$ respectively in the projective resolution (\ref{c}).
We have $$\begin{array}{ccl}\dim U_n&=&n+1\\
&=&\dim\Hom_R(P_n/(\rad P_n),\kk)\\
&=&\dim\Hom_R(P_n,\kk)\\
&=&\dim\Ext_R^n(\kk,\kk),\end{array}$$  where the last equation
follows from equation (\ref{equ mini}). So we have
$\Ext_R^*(\kk,\kk)=U$, which completes the proof of the proposition.
\qed

\subsection{Case $N\le 3$}\label{sec} In this case, the Nichols algebra
$R=\mathcal {B}(V)$ is the algebra generated by $x_1$ and $x_2$
subject to the relations
$$x_1^2x_2-(q_{_{12}}+q_{_{12}}q_{_{11}})x_1x_2x_1+q_{_{12}}^2q_{_{22}}x_2x_1^2=0,$$
$$x_2^2x_1-(q_{_{21}}+q_{_{21}}q_{_{22}})x_2x_1x_2+q_{_{21}}^2q_{_{22}}x_2x_1^2=0,$$
$$x_1^N=x_2^N=(x_1x_2-q_{_{12}}x_2x_1)^N=0.$$ The dimension of $R$ is
$N^3$.

In the rest of the paper, we set $y=x_1x_2-q_{_{12}}x_2x_1$.  From
the above relations, we obtain that
$$q_{_{21}}x_1y-yx_1=0,\;\;x_2y-q_{_{21}}yx_2=0.$$

Let $\al_1$ and $\al_2$ be the two simple roots. The element
$\al_1\al_2\al_1$ is  a reduced decomposition of the  longest
element in the Weyl group $\mc{W}$ and $\{\al_1, \al_1+\al_2,
\al_2\}$ are the positive roots. The corresponding root vectors are
just $x_1$, $y$ and $x_2$. So the set
$$\{x_1^{a_1}y^{a_2}x_2^{a_3}, \;0\se a_i<N, i=1,2,3\}$$ forms a PBW
basis of $R$. The graded algebra $\mb{G}\rr R$ corresponding to $R$
is isomorphic to the algebra generated by $x_1, y$ and $x_2$ subject
to the relations
$$x_1y=q_{_{21}}^{-1}yx_1,\;\;x_1x_2=q_{_{12}}x_2x_1,\;\;yx_2=q_{_{21}}^{-1}x_2y,$$
$$x_1^N=y^N=x_2^N=0.$$

We first show that the algebra $\Ext_R^*(\kk,\kk)$ is generated in
degree 1 and 2.

A connected graded algebra is called a \textit{$\mc{K}_2$ algebra
}\index{$\mc{K}_2$ algebra} if the algebra $\EXT^*_R(\kk,\kk)$ is
generated by $\EXT^1_R(\kk,\kk)$ and $\EXT^2_R(\kk,\kk)$. Here
$\EXT^*_R(-,-)$ denotes the functor on graded category
\cite[Definition. 1.1]{cs}.

\begin{rem}
For a finite dimensional connected algebra $R$, we have
$$\Ext_R^*(\kk,\kk)\cong \EXT_R^*(\kk,\kk).$$ In the following we just
identify them.
\end{rem}

Let $S$ be the subalgebra of $R$ generated by $x_1$ and $y$. To be
more precise, it is isomorphic to the algebra generated by $x_1$ and
$y$ subject to the relations
$$yx_1=q_{_{21}}x_1y,\;\;\;\;\;\;x_1^N=y^N=0.$$

\begin{lem}\label{k2}
The algebra $R=\mc{B}(V)$ is $\mathcal {K}_2$.
\end{lem}
\proof  The algebra $R$ is isomorphic to the graded Ore extension
$R\cong S[x_2;\sigma,\delta]$, where $\sigma$ is the graded algebra
automorphism of $S$ defined by $\sigma(x_1)=q_{_{12}}^{-1}x_1$ and
$\sigma(y)=q_{_{21}}y$ and $\delta$ is the  degree $+1$ graded
$\sigma$-derivation of $S$ defined by $\dt(x_1)=-q_{_{12}}^{-1}y$
and $\dt(y)=0$. By \cite[Thm 10.2]{cs}, the $\mc{K}_2$ property is
preserved under graded Ore extension. From \cite[Thm. 4.1]{mpsw}, we
can see that $S$ is $\mc{K}_2$. Therefore, the algebra $R$ is
$\mc{K}_2$.\qed

 The subalgebra $S$ is  a normal subalgebra of $R$ (we refer to
\cite[Appendix]{gk} for the definition of normal subalgebras). Now
set $\overline{R}=R/(RS^+)$, where $S^+$ is the augmentation ideal
of $S$. That is, $\overline{R}=k[x_2]/(x_2^{N})$.  We use the
Hochschild-Serre spectral sequence (cf. \cite{gk})
\begin{equation}E_2^{pq}=\Ext_{\overline{R}}^p(\kk,\Ext_S^q(\kk,\kk))\Longrightarrow\Ext_R^{p+q}(\kk,\kk)\end{equation} to
calculate the $\Ext$ algebra of $R$. We show that $E_2=E_{\infty}$.

The spectral sequence is constructed as follows. Let
$$\cdots\ra Q_{_1}\ra Q_0\ra
\kk\ra 0$$ and $$\cdots\ra P_1\ra P_0\ra \kk\ra 0$$ be free
 resolutions of $_{\overline{R}}\kk$ and $_R\kk$ respectively.   There is a natural $\overline{R}$-module action
on $\Hom_S(P_q,\kk)$ for $q\le 0$. We form a double complex
$$E_0^{pq}=\Hom_{\overline{R}}(Q_p,\Hom_S(P_q,\kk)).$$ By taking the
vertical homology and then the horizontal homology, we have
$$E_1^{pq}=\Hom_{\overline{R}}(Q_p,\Ext_S^q(\kk,\kk))$$ and
$$E_2^{pq}=\Ext_{\overline{R}}^p(\kk,\Ext_S^q(\kk,\kk)).$$

Now we construct a free resolution of $\kk$ over $R$, which is a
filtered complex. The corresponding graded complex is the minimal
projective resolution of $\kk$ over $\mb{G}\rr R$.

Let $\sigma, \tau:\mathbb{N}\ra\mathbb{N}$ be the functions defined
by

$$\sigma(a)=\begin{cases}1,& \text {if $a$ is odd};\\
N-1, &\text {if $a$ is even}
\end{cases}$$and
$$\tau(a)=\begin{cases}\frac{a-1}{2}N+1,& \text {if $a$ is odd};\\
\frac{a}{2}N, &\text {if $a$ is even}.
\end{cases}$$

Let \begin{equation}\label{proj}P_{\bullet}:\cdots\ra
P_n\xra{\partial_n} P_{n-1}\cdots\ra P_1\ra P_0\end{equation} be a
complex of free $R$-modules  constructed as follows. For each triple
$(a_1,a_2,a_3)$, let $\Phi(a_1,a_2,a_3)$ be a free generator for
$P_n$ with $n=a_1+a_2+a_3$. Set {\small$$P_n=\op_{a_1+a_2+a_3=n}
R\Phi(a_1,a_2,a_3)(\tau(a_1)\!+2\tau(a_2)\!+\tau(a_3), \tau(a_3),\tau(a_2),\tau(a_1)).$$} Here, (-,-,-,-) denotes the degree shift. The differentials are defined by {\small $$\partial(\Phi(a_1,a_2, a_3))=\begin{cases}(\dt_1+\dt_2+\dt_3)(\Phi(a_1,a_2, a_3)), &\text{if $a_2$ is odd};\\
(\dt_1+\dt_2+\tilde{\dt}_2+\dt_3)(\Phi(a_1,a_2, a_3)),& \text{if
$a_2$ is even}.\end{cases}$$} The maps $\dt_i$, $1\se i\se 3$ and
$\tilde{\dt}_2$ are defined as follows.

Put {\small$$\begin{array}{ccl}\dt_1(\Phi(a_1,a_2,a_3))&=&x_1^{\sigma(a_1)}\Phi(a_1-1,a_2,a_3),\;\;\; \t{ if } a_1>0;\\\dt_2(\Phi(a_1,a_2,a_3))&=&(-1)^{a_1}q_{_{21}}^{-\sigma(a_2)\tau(a_1)}y^{\sigma(a_2)}\Phi(a_1,a_2-1,a_3), \;\;\;\t{ if }a_2>0;\\
\dt_3(\Phi(a_1,a_2,a_3))&=&(-1)^{a_1+a_2}q_{_{12}}^{\sigma(a_3)\tau(a_1)}q_{_{21}}^{-\sigma(a_3)\tau(a_2)}x_2^{\sigma(a_3)}\Phi(a_1,a_2,a_3-1),\t{ if }a_3>0;\\
\tilde{\dt}_2(\Phi(a_1,a_2,a_3))&=&D\Phi(a_1-1,a_2+1,a_3-1),\;\;\;\t{
if }a_1, a_3>0, \text{$a_2$ is even},\end{array}$$} where $D$ is an
element in $R$ such that
$$Dy=-q_{_{21}}^{\tau(a_1-1)}q_{_{12}}^{\sigma(a_3)\tau(a_1-1)}q_{_{21}}^{-\sigma(a_3)\tau(a_2)}[x_1^{\sigma(a_1)},x_2^{\sigma(a_3)}]_c.$$
The existence of such element $D$ will be  explained in Lemma
\ref{d2}. For $i=1,2,3$, if  $a_i=0$, set
$\dt_i(\Phi(a_1,a_2,a_3))=0$. If $a_1=0$ or $a_3=0$, set
$\tilde{\dt}_2(\Phi(a_1,a_2,a_3))=0$.

\begin{lem}\label{d2} The element $y$ is a right divisor of
$[x_1^{\sigma(a_1)},x_2^{\sigma(a_3)}]$.
\begin{enumerate}
\item[(1)] If $a_1, a_3>0$ are odd, then
$$ \tilde{\dt}_2(\Phi(a_1,a_2,a_3))=-q_{_{21}}^{-\frac{a_2}{2}N}\Phi(a_1-1,a_2+1,a_3-1).$$

\item[(2)]  If  $a_1>0$ is odd and $a_3>0$ is even, then
{\small $$\begin{array}{cl}&\tilde{\dt}_2(\Phi(a_1,a_2,a_3))\\
=&q_{_{12}}^{(N-1)\frac{a_1-1}{2}N}q_{_{21}}^{-(N-1)\frac{a_2}{2}N}\bar{q}q_{_{21}}^{-(N-2)}q_{_{21}}^{\frac{a_1-1}{2}N}x_2^{N-2}\Phi(a_1-1,a_2+1,a_3-1).\end{array}$$}

\item[(3)] If $a_1>0$ is even and $a_3>0$ is odd, then
$$ \tilde{\dt}_2(\Phi(a_1,a_2,a_3))
=q_{_{21}}^{-\frac{a_2}{2}N}x_1^{N-2}\Phi(a_1-1,a_2+1,a_3-1).$$

\item[(4)]  If $a_1, a_3>0$ are even, then
{\small $$\begin{array}{ccl}&&\tilde{\dt}_2(\Phi(a_1,a_2,a_3))\\&=&-q_{_{12}}^{(N-1)(\frac{a_1-2}{2}N+1)}q_{_{21}}^{-(N-1)\frac{a_2}{2}N}q_{_{21}}^{\frac{a_1-2}{2}N+1}\\&&\hspace{0.4cm}(k_1x_1^{N-2}x_2^{N-2}+\cdots+k_{N-2}y^{N-3}x_1x_2+k_{N-1}y^{N-2})\Phi(a_1-1,a_2-1,a_3-1)\\
&=&-q_{_{12}}^{(N-1)(\frac{a_1-2}{2}N+1)}q_{_{21}}^{-(N-1)\frac{a_2}{2}N}q_{_{21}}^{\frac{a_1-2}{2}N+1}\\&&\hspace{0.4cm}(l_1x_2^{N-2}x_1^{N-2}+\cdots+l_{N-2}y^{N-3}x_2x_1+l_{N-1}y^{N-2})\Phi(a_1-1,a_2-1,a_3-1),\end{array}$$}
where
$$\begin{array}{ccl}[x_1^{N-1},x_2^{N-1}]_c&=&k_1yx_1^{N-2}x_2^{N-2}+\cdots+k_{N-2}y^{N-2}x_1x_2+k_{N-1}y^{N-1}\\
&=&l_1yx_2^{N-2}x_1^{N-2}+\cdots+l_{N-2}y^{N-2}x_2x_1+l_{N-1}y^{N-1},\end{array}$$
with $k_i,l_i\in \kk, 1\se i\se N-1$.
\end{enumerate}
\end{lem}
\proof (1) is easy to see. (2) and (3) follow from the following two
equations,
$$[x_1^{N-1},x_2]_c=(1+\bar{q}^{-1}+\cdots+\bar{q}^{-N+1})x_1^{N-2}y=-\bar{q}x_1^{N-2}y$$ and
$$[x_1,x_2^{N-1}]_c=(1+\bar{q}^{-1}+\cdots+\bar{q}^{-N+1})yx_2^{N-2}=-\bar{q}yx_2^{N-2}=-\bar{q}q_{_{21}}^{2-N}x_2^{N-2}y.$$
For (4), by Lemma \ref{bas} below, both
$\{x_1^{a_1}y^{a_2}x_2^{a_3}\}$ and $\{x_2^{a_3}y^{a_2}x_1^{a_1}\},
\;0\se a_i<N$, $i=1,2,3$, are bases of $R$. Using an easy induction,
we can see that $[x_1^{N-1},x_2^{N-1}]_c$ can be expressed as
$$\begin{array}{ccl}[x_1^{N-1},x_2^{N-1}]_c&=&x_1^{N-1}x_2^{N-1}-q_{_{12}}^{(N-1)^2}x_2^{N-1}x_1^{N-1}\\
&=&k_1yx_1^{N-2}x_2^{N-2}+\cdots+k_{N-2}y^{N-2}x_1x_2+k_{N-1}y^{N-1}
\\
&=&l_1yx_2^{N-2}x_1^{N-2}+\cdots+l_{N-2}y^{N-2}x_2x_1+l_{N-1}y^{N-1},
\end{array}$$ with $k_i,l_i\in \kk$, $1\se i\se N-1$. Observe that $y$ commutes with $x_1^tx_2^t$ and $x_2^tx_1^t$ for $t\le 0$.  Then the result follows. \qed

\begin{lem}\label{bas}
Both the sets $$\{x_2^{a_3}y^{a_2}x_1^{a_1}\} \;\;\text{and}
\;\;\{x_1^{a_1}y^{a_2}x_2^{a_3}\},$$ $0\se a_i<N, i=1,2,3$ form
bases of the algebra $R$.
\end{lem}
\proof It is clear for the set $\{x_1^{a_1}y^{a_2}x_2^{a_3}\}$. For
the set $\{x_2^{a_3}y^{a_2}x_1^{a_1}\}$, it is easy to see that
$$
x_2^{a_3}y^{a_2}x_1^{a_1}
=q_{_{21}}^{a_1a_2+a_2a_3}q_{_{12}}^{-a_1a_3}x_1^{a_1}y^{a_2}x_2^{a_3}+\sum_{i=1}^{min\{a_1,a_2,a_3\}}k_ix_1^{a_1-i}y^{a_2-i}x_2^{a_3-i}$$
and $$ x_1^{a_1}y^{a_2}x_2^{a_3}
=q_{_{21}}^{-a_1a_2-a_2a_3}q_{_{12}}^{a_1a_3}x_2^{a_3}y^{a_2}x_1^{a_1}+\sum_{i=1}^{min\{a_1,a_2,a_3\}}l_ix_2^{a_3-i}y^{a_2-i}x_1^{a_1-i}$$
with each $k_i,l_i\in \kk$. So $\{x_2^{a_3}y^{a_2}x_1^{a_1}\}$ also
form a basis of $R$.  \qed

\begin{prop}\label{projres}
The complex (\ref{proj}) is a projective resolution of $\kk$ over
$R$, the corresponding graded complex is the minimal projective
resolution of $\kk$ over $\mb{G}\mbb{r}R$.
\end{prop}
\proof  It is routine to check that (\ref{proj}) is indeed a
complex. We see it in   Appendix \ref{exta2. append com}. The
differentials preserve the filtration and the corresponding graded
complex is just the minimal projective resolution of $\kk$ over
$\mb{G}\mbb{r}R$ as constructed in \cite[Sec. 4]{mpsw}. Since the
filtration is finite,  the complex $P_{\bullet}$ is exact by
\cite[Chapter 2, Lemma 3.13]{b}. Therefore, $P_{\bullet}$ is a free
resolution of $\kk$ over $R$.\qed

In the following, we will forget the  shifting on the modules in the
complex (\ref{proj}). It is clear that it is still a projective
resolution of $\kk$ over $R$. The only difference is that the
differentials are not of degree 0. We denote this complex by
$P_{\bullet}$ as well.

It is well-known that the following complex is the minimal
projective resolution of $\kk$ over $\overline{R}=k[x_2]/(x_2^{N})$.
$$Q_{\bullet}:\cdots\ra\overline{R}\xra{x_2^{N-1}}\overline{R}\xra{x_2}\overline{R}\xra{x_2^{N-1}}\overline{R}\xra{x_2}\overline{R}\ra \kk.$$

Therefore, we have
$$\begin{array}{ccl}E_0^{pq}&=&\Hom_{\overline{R}}(Q_p,\Hom_S(P_q,\kk))\\
&=&\Hom_S(\op_{a_1+a_2+a_3=q}R\Phi(a_1,a_2,a_3),\kk)\\
&=&\op_{a_1+a_2+a_3=q}\overline{R}\Phi(a_1,a_2,a_3),
\end{array}$$
since $\Hom_S(R,\kk)\cong \overline{R}$. The double complex reads as
 follows

$$
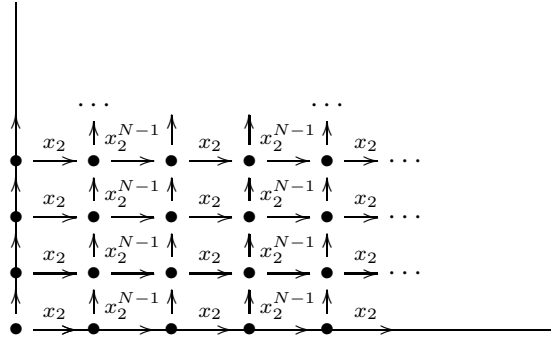
{{\xymatrix@!R@!C@=3mm{&&&&&\\
&&&&&\\
&\cdots&&&\cdots&\\
\bullet\ar[r]^{x_2}\ar[u]&\bullet\ar[r]^{x_2^{N-1}}\ar[u]&\bullet\ar[r]^{x_2}\ar[u]&\bullet\ar[r]^{x_2^{N-1}}\ar[u]&\bullet\ar[r]^{x_2}\ar[u]&\cdots\\
\bullet\ar[r]^{x_2}\ar[u]&\bullet\ar[r]^{x_2^{N-1}}\ar[u]&\bullet\ar[r]^{x_2}\ar[u]&\bullet\ar[r]^{x_2^{N-1}}\ar[u]&\bullet\ar[r]^{x_2}\ar[u]&\cdots\\
\bullet\ar[r]^{x_2}\ar[u]&\bullet\ar[r]^{x_2^{N-1}}\ar[u]&\bullet\ar[r]^{x_2}\ar[u]&\bullet\ar[r]^{x_2^{N-1}}\ar[u]&\bullet\ar[r]^{x_2}\ar[u]&\cdots\\
\bullet\ar[r]\ar[r]^{x_2}\ar@{-}[uuuuuu]\ar@{-}[rrrrrrr]\ar[u]&\bullet\ar[r]^{x_2^{N-1}}\ar[u]&\bullet\ar[r]^{x_2}\ar[u]&\bullet\ar[r]^{x_2^{N-1}}\ar[u]&\bullet\ar[r]^{x_2}\ar[u]&&&\\}}}
$$
The vertical differentials are induced from the differentials of the
complex (\ref{proj}).

By taking the vertical homology,  we have
$E_1^{pq}=\Hom_{\overline{R}}(Q_p,\Ext_S^q(\kk,\kk))$. Following
from \cite{mpsw}, the algebra $\Ext_S^*(\kk,\kk)$ is generated by
$\mk{u}_1$, $\mk{u}_y$, $\mk{w}_1$ and $\mk{w}_y$, where $\deg
\mk{u}_1=\deg \mk{u}_y=2$ and $\deg \mk{w}_1=\deg \mk{w}_y=1$,
subject to the relations
$$\mk{w}_y\mk{w}_1=-q_{_{21}}\mk{w}_1\mk{w}_y,\;\;\mk{w}_1^2=\mk{w}_y^2=0,$$
$$\mk{w}_y\mk{u}_1=q_{_{21}}^{N}\mk{u}_1\mk{w}_y,\;\;\mk{w}_1\mk{u}_1=\mk{u}_1\mk{w}_1,\;\;\mk{w}_y\mk{u}_y=\mk{u}_y\mk{w}_y,\;\;\mk{w}_1\mk{u}_y=q_{_{21}}^{-N}\mk{u}_y\mk{w}_1,$$
$$\mk{u}_y\mk{u}_1=q_{_{21}}^{N^2}\mk{u}_1\mk{u}_y.$$

We use the notations $\mk{u}_i$ and $\mk{w}_i$ in place of the
notations $\xi_i$ and $\eta_i$ used in \cite{mpsw}.  Note that
$\mk{w}_1^2=\mk{w}_y^2=0$ holds since we assume that the
characteristic of the field $\kk$ is 0. It should also be noticed
that the Ext algebra in \cite{mpsw} is the opposite algebra here.

As described in the appendix of \cite{gk}, there is an action of
$\overline{R}$ on $\Ext_S^*(\kk,\kk)$ given by
$$x_2(\mk{u}_y)=x_2(\mk{u}_1)=0, \;\;x_2(\mk{w}_y)=\mk{w}_1,\t{ and } x_2(\mk{w}_1)=0.$$
This action is a derivation on $\Ext_S^*(\kk,\kk)$. That is,
$x_2(\mk{u}\mk{w})=x_2(\mk{u})\mk{w}+\mk{u}x_2(\mk{w})$ for
$\mk{u},\mk{w}\in \Ext_S^*(\kk,\kk) $.

The following lemma gives a basis of
$\Ext_{\overline{R}}^p(\kk,\Ext_S^q(\kk,\kk))$.

\begin{lem}\label{eebasis}
As a vector space, $\Ext_{\overline{R}}^p(\kk,\Ext_S^q(\kk,\kk))$
has a basis as follows
$$\begin{cases} \mk{u}_1^i\mk{u}_y^j\mk{w}_1,   \;\:\;\;\;\;\;\;\;\;\;2(i+j)+1=q, & \text{q is odd and p is even};\\
\mk{u}_1^i\mk{u}_y^j\mk{w}_y,   \;\:\;\;\;\;\;\;\;\;\;2(i+j)+1=q, & \text{q is odd and p is odd};\\
\mk{u}_1^i\mk{u}_y^j(\mk{w}_1\mk{w}_y)^k,\;\;k=0,1 \t{ and }
2(i+j)+2k=q, & \text{q is even}.
\end{cases}$$
\end{lem}
\proof Let $E=\Ext_S^*(\kk,\kk)$. The lemma follows directly from
the following facts:

(i) If  $q$ is odd, then $\{\mk{u}_1^i\mk{u}_y^j\mk{w}_1|i,j\le
0,2(i+j)+1=q\}$ forms a basis of $x_2E^q$ and $\{e\in E^q|x_2e=0\}$.

(ii) If  $q$ is even, then $x_2E^q=0$.

(iii) $x_2^{N-1}E=0$. \qed

\begin{prop}\label{spec}
The spectral sequence
$$E_2^{p,q}=\Ext_{\overline{R}}^p(\kk,\Ext_S^q(\kk,\kk))\Longrightarrow\Ext_R^{p+q}(\kk,\kk)$$ satisfies $E_2=E_{\infty}$.
\end{prop}
\proof The elements $\mk{u}_1^i\mk{u}_y^j\mk{w}_y$ and
  $\mk{u}_1^i\mk{u}_y^j\mk{w}_1$ are represented by
$$x_2^{N-2}\Phi(2i+1,2j,0)+q_{_{12}}^{-(j+1)}x_2^{N-1}\Phi(2i,2j+1,0)$$
and  $$x_2^{N-1}\Phi(2i+1,2j,0),$$ while $\mk{u}_1^i\mk{u}_y^j$ and
$\mk{u}_1^i\mk{u}_y^j\mk{w}_1\mk{w}_y$ are represented by
$$x_2^{N-1}\Phi(2i,2j,0)\t{ and }x_2^{N-1}\Phi(2i+1,2j+1,0)$$ in $E_0$.
In other words, all the elements in $E_0$ representing  the elements
in $E_2$ are mapped to 0 under the horizontal differentials. We
conclude that $E_2=E_{\infty}$. \qed

%Now we shall analyze the structure of $E_2^{pq}=\Ext_{\overline{R}}^p(k,\Ext_S^q(\kk,\kk))$.

We now can determine the dimension of $\Ext^*_R(\kk,\kk)$. This
dimension  depends on the parity of $n$.

\begin{cor}\label{dim} We have
$$\dim \Ext^n_R(\kk,\kk)=\begin{cases} \frac{3n^2+8n+5}{8}
,& \text{if  $n$ is odd;}\\
\frac{3n^2+10n+8}{8}, &  \text{if $n$ is even.}
\end{cases}$$
\end{cor}
\proof  Set
$E^n=\op_{p+q=n}E_2^{pq}=\op_{p+q=n}\Ext_{\overline{R}}^p(\kk,\Ext_S^q(\kk,\kk))$.
By  Lemma \ref{eebasis}, we can illustrate the dimensions of
$E_2^{pq}$ with the following table:
$${\small\xymatrix@!R@!C@=3mm{
&&&&&\\
&&\cdots&\cdots&&\\
4\ar@{-}[uu]&4&4&4&4&4&\cdots& \\
7\ar@{-}[u]&7&7&7&7&7&\cdots& \\
3\ar@{-}[u]&3&3&3&3&3&\cdots& \\
5\ar@{-}[u]&5&5&5&5&5&\cdots& \\
2\ar@{-}[u]&2&2&2&2&2&\cdots& \\
3\ar@{-}[u]&3&3&3&3&3&\cdots& \\
1\ar@{-}[u]&1&1&1&1&1&\cdots& \\
1\ar@{-}[r]\ar@{-}[u]&1\ar@{-}[r]&\ar@{-}[r]1\ar@{-}[r]&\ar@{-}[r]1&\ar@{-}[r]1&\ar@{-}[r]1&\cdots&.}}
$$%
Therefore, when $n$ is odd, $$\begin{array}{ccl}\dim
E^n&=&(1+2+\cdots+\frac{n+1}{2}+1+3+\cdots+n)\\
&=&\frac{3n^2+8n+5}{8}.
\end{array}$$
When $n$ is even,
 $$\begin{array}{ccl}\dim
E^n&=&(1+2+\cdots+\frac{n}{2}+1+3+\cdots+n+1)\\
&=&\frac{3n^2+10n+8}{8}.
\end{array}$$
By Proposition \ref{spec}, we have $E_2=E_{\infty}$, so
$\dim\Ext^n_R(\kk,\kk)=\dim E^n$. This completes the proof. \qed

Now we give the first segment of the minimal projective resolution
of a Nichols algebra of type $A_2$.

The algebra $R$ is a local algebra. Thus projective $R$-modules are
free. Let $$R^{n_4}\ra R^{n_3}\ra R^{n_2}\ra R^{n_1}\ra R^{n_0}\ra
\kk\ra 0$$ be the first segment of the minimal projective
resolution. Since $\kk$ is a simple module, we have
$$\begin{array}{ccl}\dim\Ext^i_R(\kk,\kk)&=&\dim\Hom_R(R^{n_i},\kk)\\&=&\dim\Hom_R((R/(\rad
R))^{n_i},\kk)\\&=&n_i.\end{array}$$ From the computation of the
dimensions of $\Ext_R^*(\kk,\kk)$ in Corollary \ref{dim}, we can see
that the minimal projective resolution begins as $$R^{12}\ra R^7\ra
R^5\ra R^2\ra R\ra \kk\ra 0.$$ We give the differentials in the
following proposition.

As in the construction of $\tilde{\dt}_2$ in \S\ref{sec}, let
$\overline{D}$ be the element in  $R$  such that
$\overline{D}y=[x_1^{N-1},x_2^{N-1}]_c$.

\begin{prop}\label{res1}
Let $R$ be a Nichols algebra of type $A_2$. The following sequence
provides the first segment of the minimal projective resolution of
$\kk$ over $R$,
\begin{equation}\label{res}R^{12}\xra{d_4} R^7\xra{d_3} R^5\xra{d_2}
R^2\xra{d_1} R\ra \kk\ra 0,\end{equation} where the differentials
are given by the following matrices:
\end{prop}
$$d_1={\small\left(\begin{array} {c}x_1\\x_2
\end{array}\right)},$$

$$d_2={\small\left(\begin{array} {cc}x_1^{N-1}&0\\-(q_{_{12}}+\bar{q}q_{_{12}})x_1x_2+\bar{q}q_{_{12}}^2x_2x_1&x_1^2\\
-q_{_{12}}y^{N-1}x_2&y^{N-1}x_1\\
x_2^2&\bar{q}q_{_{21}}^2x_1x_2-(q_{_{21}}+\bar{q}q_{_{21}})x_2x_1\\
0&x_2^{N-1}
\end{array}\right)},$$

$$d_3={\small\left(\begin{array} {ccccc}x_1&0&0&0&0\\
q_{_{12}}^Nx_2&x_1^{N-2}&0&0&0\\
0&0&x_2&q_{_{12}}q_{_{21}}^{N-1}y^{N-1}&0\\
0&x_2&0&x_1&0\\
0&-q_{_{21}}^{1-N}y^{N-1}&x_1&0&0\\
0&0&0&q_{_{12}}^Nx_2^{N-2}&x_1\\
0&0&0&0&x_2
\end{array}\right)},$$

$$d_4=
\left(\begin{array}{c;{2pt/2pt}c}
A_1&A_2\\
\hdashline[2pt/2pt]
A_3&A_4\\
\end{array}\right),$$

where $$A_1={\tiny\left(\begin{array} {cccc}
x_1^{N-1}&0&0&0\\
-(q_{_{12}}^{1+N}+\bar{q}q_{_{12}}^{1+N})x_1x_2+\bar{q}q_{_{12}}^{2+N}x_2x_1&x_1^2&0&0\\
-q_{_{12}}^{1+N}y^{N-1}x_2&y^{N-1}x_1&0&0\\
q_{_{12}}^Nx_2^2&\bar{q}q_{_{21}}^2x_1x_2-(q_{_{21}}+\bar{q}q_{_{21}})x_2x_1&0&q_{_{21}}^Nx_1^{N-1}\\
0&x_2^{N-1}&0&-q_{_{12}}^{-N^2+2N}\overline{D}
\end{array}\right),}$$

$$A_2={\tiny\left(\begin{array} {ccc}
0&0&0\\
0&0&0\\
q_{_{12}}^{-N^2+N}x_1^{N-1}&0&0\\
0&0&0\\
0&q_{_{12}}^{-N^2+N}x_1^{N-1}&0
\end{array}\right)},$$

$$A_3={\tiny\left(\begin{array} {cccc}
0&0&x_1^2&-\bar{q}^{-1}q_{_{12}}^{N}y^{N-1}x_1\\
0&0&y^{N-1}x_1&0\\
0&0&\bar{q}q_{_{21}}^2x_1x_2-(q_{_{21}}+\bar{q}q_{_{21}})x_2x_1&q_{_{21}}^{N-1}y^{N-1}x_2\\
0&0&0&q_{_{12}}^{2N}x_2^{N-1}\\
0&0&q_{_{12}}^{N^2}x_2^{N-1}&0\\
0&0&0&0\\0&0&0&0
\end{array}\right)},$$

{\small$$A_4={\tiny\left(\begin{array} {ccc}
-(q_{_{12}}+\bar{q}q_{_{12}})x_1x_2+\bar{q}q_{_{12}}^2x_2x_1&0&0\\
-q_{_{12}}y^{N-1}x_2&0&0\\
x_2^2&0&0\\
0&-(q_{_{12}}+\bar{q}q_{_{12}})x_1x_2+\bar{q}q_{_{12}}^2x_2x_1&x_1^2\\
0&-q_{_{12}}y^{N-1}x_2&y^{N-1}x_1\\
0&x_2^2&\bar{q}q_{_{21}}^2x_1x_2-(q_{_{21}}+\bar{q}q_{_{21}})x_2x_1\\
0&0&x_2^{N-1}
\end{array}\right)}.$$}
\proof It is routine to check that (\ref{res}) is indeed a complex.
But we need to mention that the following two equations hold
$$\overline{D}x_1-x_1^{N-1}x_2^{N-2}=0,$$
$$x_2^{N-1}x_1^{N-2}-q_{_{12}}^{-N^2+2N}\overline{D}x_2=0.$$ These equations follow
from Lemma \ref{bas} and the equations
$$[x_1^{N-1},x_2^{N-1}]_cx_1=yx_1^{N-1}x_2^{N-2},$$
$$[x_1^{N-1},x_2^{N-1}]_cx_2=q_{_{12}}^{N^2-2N}yx_2^{N-1}x_1^{N-2}.$$

The complex (\ref{res}) is homotopically equivalent to the first
segment of the resolution $P_\bullet$ (without shifting) constructed
in Section 2.  Therefore, it is exact. \qed

\begin{rem}
In \cite[Theorem 6.1.3]{mw}, the authors give a set of linearly
independent 2-cocycles  on $R$, indexed by the positive roots. In
the  resolution (\ref{res}), the functions dual to $(1,0,0,0,0)$,
$(0,0,1,0,0)$ and $(0,0,0,0,1)$ are just those 2-cocycles,
corresponding to the positive roots $\al_1$, $\al_1+\al_2$ and
$\al_2$ respectively.
\end{rem}

Now we give our main theorems about the structure of the $\Ext$
algebra of a Nichols algebra of type $A_2$.
\begin{thm}\label{1}
Let $R$ be a Nichols algebra of type $A_2$ with  $N=3$, then
$\Ext_R^*(\kk,\kk)$ is generated by $\mk{a}_i$, $\mk{b}_i$,
$\mk{c}_i$, $i=1,2$ and $\mk{b}_y$ with
$$\deg \mk{a}_i=1,\;\;\deg \mk{b}_i=\deg \mk{b}_y=\deg \mk{c}_i=2,$$ subject to
the relations
$$\mk{a}_1^2=\mk{a}_2^2=0,\;\;\;\mk{a}_1\mk{a}_2=\mk{a}_2\mk{a}_1=0,$$
$$\mk{a}_1\mk{b}_1=\mk{b}_1\mk{a}_1,\;\;\mk{a}_1\mk{b}_y=q_{_{12}}^{3}\mk{b}_y\mk{a}_1,\;\;\mk{a}_1\mk{b}_2=q_{_{12}}^3\mk{b}_2\mk{a}_1,$$
$$\mk{a}_1\mk{c}_1=\bar{q}^2q_{_{12}}\mk{c}_1\mk{a}_1,\;\;\mk{a}_1\mk{c}_2=\bar{q}q_{_{12}}^2\mk{c}_2\mk{a}_1,$$
$$q_{_{12}}^3\mk{a}_2\mk{b}_1=\mk{b}_1\mk{a}_2,\;\;q_{_{12}}^{3}\mk{a}_2\mk{b}_y=\mk{b}_y\mk{a}_2,\;\;\mk{a}_2\mk{b}_2=\mk{b}_2\mk{a}_2,$$
$$\mk{a}_2\mk{c}_1=\bar{q}q_{_{21}}^2\mk{c}_1\mk{a}_2,\;\;\bar{q}^2q_{_{12}}\mk{a}_2\mk{c}_2=\mk{c}_2\mk{a}_2,$$
$$\bar{q}^2q_{_{12}}\mk{a}_2\mk{b}_1=\mk{a}_1\mk{c}_1,\;\;\mk{a}_1\mk{b}_2=\bar{q}^2q_{_{12}}\mk{a}_2\mk{c}_2,\;\;\mk{c}_1\mk{a}_2=\mk{c}_2\mk{a}_1,$$
$$\mk{b}_1\mk{c}_2=q_{_{12}}^6\mk{c}_1^2,\;\;q_{_{12}}^{6}\mk{b}_2\mk{c}_1=\mk{c}_2^2,\;\;\mk{b}_1\mk{b}_2=q_{_{12}}^3\mk{c}_1\mk{c}_2,\;\;\mk{c}_1\mk{c}_2=q_{_{12}}^3\mk{c}_2\mk{c}_1,$$
$$\mk{b}_1\mk{b}_y=q_{_{12}}^9\mk{b}_y\mk{b}_1,\;\;\mk{b}_1\mk{b}_2=q_{_{12}}^9\mk{b}_2\mk{b}_1,\;\;\mk{b}_y\mk{b}_2=q_{_{12}}^9\mk{b}_2\mk{b}_y,$$
$$q_{_{12}}^3\mk{c}_1\mk{b}_1=\mk{b}_1\mk{c}_1,\;\;\mk{c}_1\mk{b}_y=q_{_{12}}^3\mk{b}_y\mk{c}_1,\;\;\mk{c}_1\mk{b}_2=q_{_{12}}^6\mk{b}_2\mk{c}_1,$$
$$q_{_{12}}^6\mk{c}_2\mk{b}_1=\mk{b}_1\mk{c}_2,\;\;q_{_{12}}^3\mk{c}_2\mk{b}_y=\mk{b}_y\mk{c}_2,\;\;\mk{c}_2\mk{b}_2=q_{_{12}}^3\mk{b}_2\mk{c}_2.$$
\end{thm}

\begin{thm}\label{2}
Let $R$ be a Nichols algebra of type $A_2$ with  $N>3$, then
$\Ext_R^*(\kk,\kk)$ is generated by $\mk{a}_i$, $\mk{b}_i$ and
$\mk{c}_i$, $i=1,2$  and $\mk{b}_y$ with
$$\deg \mk{a}_i=1,\;\;\deg \mk{b}_i=\deg \mk{b}_y=\deg \mk{c}_i=2,$$  subject to
the relations
$$\mk{a}_1^2=\mk{a}_2^2=0,\;\;\;\mk{a}_1\mk{a}_2=\mk{a}_2\mk{a}_1=0,$$
$$\mk{a}_1\mk{b}_1=\mk{b}_1\mk{a}_1,\;\;\mk{a}_1\mk{b}_y=q_{_{12}}^N\mk{b}_y\mk{a}_1,\;\;\mk{a}_1\mk{b}_2=q_{_{12}}^N\mk{b}_2\mk{a}_1,$$
$$q_{_{12}}^N\mk{a}_2\mk{b}_1=\mk{b}_1\mk{a}_2,\;\;q_{_{12}}^{N}\mk{a}_2\mk{b}_y=\mk{b}_y\mk{a}_2,\;\;\mk{a}_2\mk{b}_2=\mk{b}_2\mk{a}_2,$$
$$\mk{a}_1\mk{c}_2=\bar{q}q_{_{12}}^2\mk{c}_2\mk{a}_1,\;\;\mk{a}_2\mk{c}_1=\bar{q}q_{_{21}}^2\mk{c}_1\mk{a}_2,$$
$$\mk{a}_1\mk{c}_1=\mk{c}_1\mk{a}_1=\mk{c}_2\mk{a}_2=\mk{a}_2\mk{c}_2=0,\;\;\mk{c}_1\mk{a}_2=\mk{c}_2\mk{a}_1,$$
$$\mk{c}_1^2=\mk{c}_2^2=\mk{c}_1\mk{c}_2=\mk{c}_2\mk{c}_1=0,$$
$$\mk{b}_1\mk{b}_y=q_{_{12}}^{N^2}\mk{b}_y\mk{b}_1,\;\;\mk{b}_1\mk{b}_2=q_{_{12}}^{N^2}\mk{b}_2\mk{b}_1,\;\;\mk{b}_y\mk{b}_2=q_{_{12}}^{N^2}\mk{b}_2\mk{b}_y,$$
$$q_{_{12}}^N\mk{c}_1\mk{b}_1=\mk{b}_1\mk{c}_1,\;\;\mk{c}_1\mk{b}_y=q_{_{12}}^N\mk{b}_y\mk{c}_1,\;\;\mk{c}_1\mk{b}_2=q_{_{12}}^{2N}\mk{b}_2\mk{c}_1,$$
$$q_{_{12}}^{2N}\mk{c}_2\mk{b}_1=\mk{b}_1\mk{c}_2,\;\;q_{_{12}}^N\mk{c}_2\mk{b}_y=\mk{b}_y\mk{c}_2,\;\;\mk{c}_2\mk{b}_2=q_{_{12}}^N\mk{b}_2\mk{c}_2.$$
\end{thm}

\it{Proof of Theorems \ref{1} and \ref{2}} We prove Theorem \ref{1}.
Theorem \ref{2} can be proved similarly.  Consider the minimal
resolution (\ref{res}) showed in Proposition \ref{res1}, we have
$\Ext_R^1(\kk,\kk)=\Hom_R(R^2,\kk)$ and
$\Ext_R^2(\kk,\kk)=\Hom_R(R^5,\kk)$, since $\kk$ is a simple module.
Let $\mk{a}_1,\mk{a}_2\in \Ext_R^1(\kk,\kk)$ be the functions dual
to $(1,0)$ and $(0,1)$ respectively. Let
$\mk{b}_1,\mk{c}_1,\mk{b}_y,\mk{c}_2,\mk{b}_2\in \Ext_R^2(\kk,\kk)$
be the functions dual to $(1,0,0,0,0)$,$\cdots$, $(0,0,0,0,1)$
respectively. The relations listed in the theorem can be verified by
constructing suitable commutative diagrams, we do this in Appendix
\ref{exta2. append relation}. Let $U$ be an algebra generated by
$\mk{b}_1$, $\mk{b}_y$, $\mk{b}_2$ and $\mk{a}_i$, $\mk{c}_i$,
$i=1,2$, subject to the relations listed in the theorem. Then any
element in $U$  can be written as a linear combination of elements
of the form
$\mk{b}_1^{b_1}\mk{b}_y^{b_2}\mk{b}_2^{b_3}\mk{a}_i^{a_i}$,
$\mk{b}_1^{b_1}\mk{b}_y^{b_y}\mk{b}_2^{b_2}\mk{c}_i^{c_i}$ and
$\mk{b}_1^{b_1}\mk{b}_y^{b_y}\mk{b}_2^{b_2}\mk{c}_1\mk{a}_2$, with
$b_1$, $b_2$, $b_3\le 0$, $a_i$, $c_i\in \{0,1\}$, $i=1,2$.

 By Lemma \ref{k2}, the algebra $U$ is  a quotient of $\Ext^*_R(\kk,\kk)$.  When $n$ is
odd, $$\begin{array}{ccl}\dim U_n&=&(\frac{n-1}{2}+2)(\frac{n-1}{2}+1)+\frac{1}{2}(\frac{n-1}{2})(\frac{n-1}{2}+1)\\
&=&\frac{3n^2+8n+5}{8}. \end{array}$$ When $n$ is even,
$$\begin{array}{ccl}\dim
U_n&=&(\frac{n}{2})(\frac{n}{2}+1)+\frac{1}{2}(\frac{n}{2}+1)(\frac{n}{2}+2)\\
&=&\frac{3n^2+10n+8}{8}. \end{array}$$ It follows from  Corollary
\ref{dim}  that $\dim U_n=\dim\Ext_R^n(\kk,\kk)$, for all $n\le 0$,
so $U=\Ext^*_R(\kk,\kk)$, which completes the proof of the theorem.
\qed

\begin{rem}
In \cite[Thm 5.4]{mpsw}, the authors showed that the Ext algebra of
a Nichols algebra of finite Cartan type is braided commutative. This
coincides with the results we obtain in Theorems \ref{1} and
\ref{2}.
\end{rem}

Now we can answer the question whether the Ext algebra of a Nichols
algebra is still a Nichols algebra. In general, the answer is
negative.

\begin{prop}
The Ext algebra of a Nichols algebra of type $A_2$ is not a Nichols
algebra.
\end{prop}
\proof We consider the case $N=2$ first. Denote the Ext algebra by
$E$.  From Proposition \ref{N=2}, $E$ is generated by $\mk{a}_1$,
$\mk{a}_2$ and $\mk{b}$  subject to the relations
$$\mk{a}_2\mk{a}_1=\mk{a}_1\mk{a}_2=0,\;\;\mk{a}_1\mk{b}=\mk{b}\mk{a}_1,\;\;\mk{a}_2\mk{b}=\mk{b}\mk{a}_2.$$
If $E$ is a Nichols algebra with respect to some braided vector
space $V$, then $\mk{a}_1$, $\mk{a}_2$ and $\mk{b}$ should form a
basis of $V$. This is because as an algebra, a Nichols algebra
$\mc{B}(V)$ is generated by elements in $V$.  With relation
$\mk{a}_2\mk{a}_1=\mk{a}_1\mk{a}_2$, $\mk{a}_1\mk{b}=\mk{b}\mk{a}_1$
and $\mk{a}_2\mk{b}=\mk{b}\mk{a}_2$, the vector space $V$ is of
diagonal type. This contradicts to the relation
$\mk{a}_2\mk{a}_1=\mk{a}_1\mk{a}_2=0$. Therefore, $E$ is not a
Nichols algebra. By a similar argument, we can conclude that when
$N\le 3$, the Ext algebra is not a Nichols algebra either. \qed

However, we have the following positive result.

\begin{prop}
Let $R$ be a Nichols algebra of type $A_2$ with $N> 3$. Then
$\Ext_R^*(\kk,\kk)/\mc{N}$ is a Nichols algebra of diagonal type,
where $\mc{N}$ is the ideal generated by nilpotent elements.
\end{prop}
\proof From the proof of   Theorem \ref{2}, the elements
$\mk{b}_1^{b_1}\mk{b}_y^{b_2}\mk{b}_2^{b_3}\mk{a}_i^{a_i}$,
$\mk{b}_1^{b_1}\mk{b}_y^{b_y}\mk{b}_2^{b_2}\mk{c}_i^{c_i}$ and
$\mk{b}_1^{b_1}\mk{b}_y^{b_y}\mk{b}_2^{b_2}\mk{c}_1\mk{a}_2$, with
$b_1$, $b_2$, $b_3\le 0$, $a_i$, $c_i\in \{0,1\}$, $i=1,2$ form a
basis of $\Ext_R^*(\kk,\kk)$. With the relation listed in that
theorem, the elements
$\mk{b}_1^{b_1}\mk{b}_y^{b_2}\mk{b}_2^{b_3}\mk{a}_i $,
$\mk{b}_1^{b_1}\mk{b}_y^{b_y}\mk{b}_2^{b_2}\mk{c}_i $ and
$\mk{b}_1^{b_1}\mk{b}_y^{b_y}\mk{b}_2^{b_2}\mk{c}_1\mk{a}_2$ are
nilpotent. However, linear combination of elements
$\mk{b}_1^{b_1}\mk{b}_y^{b_2}\mk{b}_2^{b_3}$  are not nilpotent.
Then the algebra $\Ext_R^*(\kk,\kk)/\mc{N}$ is generated by
$\mk{b}_1$, $\mk{b}_y$ and $\mk{b}_2$  subject to the relations
$$\mk{b}_1\mk{b}_y=q_{_{12}}^{N^2}\mk{b}_y\mk{b}_1,\;\;\mk{b}_1\mk{b}_2=q_{_{12}}^{N^2}\mk{b}_2\mk{b}_1,\;\;\mk{b}_y\mk{b}_2=q_{_{12}}^{N^2}\mk{b}_2\mk{b}_y.$$
It is obvious that it is a Nichols algebra of diagonal type with
Cartan matrix of type $A_1\times A_1\times A_1$.\qed

The following corollary is a direct consequence from Theorem \ref{1}
and \ref{2}.

\begin{cor}\label{var}
Let $A=u(\mathcal{D}, 0, \mu)$ be a pointed Hopf algebra of type
$A_2$ with $N\le 3$ and $R=\mathcal{B}(V)$   the corresponding
Nichols algebra. Then
$$\cx_R(\kk)=\cx_A(\kk)=3.$$ In addition,   $\mathcal{V}_A(\kk)\cong \mathcal{V}_{(\mb{G}\rr R)\#\kk\bgm}(\kk)$.
\end{cor}
\proof For the  Nichols algebra $R$, the complexity
$$\cx_R(\kk)=\gm(\Ext_R^*(\kk,\kk))=3$$ follows directly from
Proposition \ref{dim} or Theorems \ref{1} and \ref{2}. By
\cite[Lemma 6.1]{mpsw}, we have $$\H^*(u(\mathcal{D}, 0,
\mu),\kk)\cong \H^*(u(\mathcal{D}, 0, 0),\kk).$$ In addition, we
also have
$$\Ext^*_{u(\mathcal{D}, 0, 0)}(\kk,\kk)\cong \Ext^*_R(\kk,\kk)^\bgm.$$
Observe that for each positive root $\al$, some power of
$\mk{b}_{\al}$ is invariant under the group action. Indeed, from the
discussion in Section 6 in \cite{mw}, each $\mk{b}_{\al}$ (denoted
by $f_\al$ there) can be expressed as a function $R^+\times R^+\ra
\kk$. Then we see that $\mk{b}_{\al}^{M_\al}$ is $\bgm$-invariant,
where $M_\al$ is the integer such that $\chi_\al^{M_\al}=\vps$.
Hence, $\gm(\H^*(u(\mathcal{D}, 0, 0),\kk)=3$, which implies that
$\cx_A(\kk)=3$. With the relations in Theorems \ref{1} and \ref{2},
we see that $$\mathcal{V}_A(\kk)\cong
\t{MaxSpec}(\kk[\mk{b}_1^{m_1},\mk{b}_y^{m_y},\mk{b}_2^{m_2}]),$$
where $m_1$, $m_y$ and $m_2$ are the least integers such that
$\mk{b}_1^{m_1}$, $\mk{b}_y^{m_y}$, $\mk{b}_2^{m_2}\in
\H^*(u(\mathcal{D}, 0, 0),\kk)$.  That is, $\mathcal{V}_A(\kk)$ is
isomorphic to the maximal spectrum of the polynomial algebra
$\kk[y_1,y_2,y_3].$ By \cite[Thm. 4.1]{mpsw} $\mathcal{V}_{\mb{G}\rr
R\#\kk\bgm}(\kk)$ is also isomorphic to the maximal spectrum of
$\kk[y_1,y_2,y_3].$ So $\mathcal{V}_A(\kk)\cong
\mathcal{V}_{\mb{G}\rr R\#\kk\bgm}(\kk)$. \qed
%\begin{rem}
%The above corollary stated that $\cx_R(k)=\cx_{R}(k)$ for Nichols algebra of type $A_2$. We conjecture that this is true for Nichols algebra of %finie Cartan type.
%\end{rem}

To end this section, we give an easy application of the main
theorems. We show that a large class of  finite dimensional pointed
Hopf algebras of finite Cartan type  are wild.

\begin{prop}\label{repre}
Let $A=u(\mathcal{D}, \lmd, \mu)$ be a pointed Hopf algebra such
that the components of the Dynkin diagram  are of
 type $A$, $D$, or $E$, except for $A_1$ and
$A_1\times A_1$, and the order $N_J>2$ for at least one component.
Then $A$ is wild.
\end{prop}
\proof In view of \cite[Thm. 3.1]{fw}, we only need to prove that
$\cx_A(\kk)\le 3$. Using \cite[Lemma 6.1]{mpsw} again, we have
$\cx_A(\kk)=\cx_{u(\mathcal{D},\lambda,0)}(\kk)$. However,
$u(\mathcal{D},\lambda,0)$ contains a Hopf subalgebra $B$ which is
of type $A_2$ with the order $N\le3$. Thus
$\cx_{u(\mathcal{D},\lambda,0)}(\kk)\le\cx_B(\kk)\le 3$ by
\cite[Prop 2.1]{fw}. \qed

We conjecture that the isomorphism $\mathcal{V}_A\cong
\mathcal{V}_{\mathbb{G}\rr R\# \kk \bgm}$ in Corollary  \ref{var}
holds for general finite dimensional pointed Hopf algebra
$A=u(\mathcal{D},\lmd, \mu)$ of finite Cartan type.

\section{Appendix}\label{exta2. append}

\subsection{}\label{exta2. append com}In this subsection, we verify that the complex (\ref{proj}) in \S
\ref{sec} is indeed a complex.

The following equations follow directly from Lemma \ref{d2}.
\begin{equation}\label{dy}Dy=\begin{cases}yD,&\text{if $a_1,a_3$ are both even or both  odd};\\
q_{_{21}}^{-N+2}yD, & \text{if $a_1$ even and $a_3$ is odd};\\
q_{_{21}}^{N-2}yD, & \text{if $a_1$ odd and $a_3$ is even}.\\
\end{cases}\end{equation}

It is clear that $\dt_i^2=0$ for $i=1,2,3$. So if $a_2$ is odd,
{\small $$\partial^2(\Phi(a_1,a_2,
a_3))=((\dt_3\dt_1+\dt_1\dt_3+\tilde{\dt}_2\dt_2)+(\dt_2\dt_3+\dt_3\dt_2)+(\dt_1\dt_2+\dt_2\dt_1))\Phi(a_1,a_2,
a_3).$$} Put  $${\small
\begin{array}{ccl}A&=&(\dt_3\dt_1+\dt_1\dt_3+\tilde{\dt}_2\dt_2)\Phi(a_1,a_2,
a_3),\\
B&=&(\dt_2\dt_3+\dt_3\dt_2)\Phi(a_1,a_2, a_3),\\
C&=&(\dt_1\dt_2+\dt_2\dt_1)\Phi(a_1,a_2, a_3).\end{array}}$$ We show
that $A=B=C=0$. {\small
 $$\begin{array}{ccl}A&=&(\dt_3\dt_1+\dt_1\dt_3+\tilde{\dt}_2\dt_2)\Phi(a_1,a_2, a_3)\\
&=&((-1)^{a_1-1+a_2}q_{_{12}}^{\sigma(a_3)\tau(a_1-1)}q_{_{21}}^{-\sigma(a_3)\tau(a_2)}[x_1^{\sigma(a_1)},x_2^{\sigma(a_3)}]_c\\
&&\hspace{5mm}+(-1)^{a_1}q_{_{21}}^{-\tau(a_1)}yD)\Phi(a_1-1,a_2,a_3-1),
\end{array}$$} where $D$ satisfies that $$Dy=-q_{_{21}}^{\tau(a_1-1)}q_{_{12}}^{\sigma(a_3)\tau(a_1-1)}q_{_{21}}^{-\sigma(a_3)\tau(a_2-1)}[x_1^{\sigma(a_1)},x_2^{\sigma(a_3)}]_c.$$ That is,
$$q_{_{12}}^{\sigma(a_3)\tau(a_1-1)}q_{_{21}}^{-\sigma(a_3)\tau(a_2-1)}[x_1^{\sigma(a_1)},x_2^{\sigma(a_3)}]_c+q_{_{21}}^{-\tau(a_1-1)}Dy=0.$$
Hence,
$$q_{_{12}}^{\sigma(a_3)\tau(a_1-1)}q_{_{21}}^{-\sigma(a_3)\tau(a_2)}[x_1^{\sigma(a_1)},x_2^{\sigma(a_3)}]_c+q_{_{21}}^{-\sigma(a_3)}q_{_{21}}^{-\tau(a_1-1)}Dy=0.$$
By equation (\ref{dy}), we have
$q_{_{21}}^{-\sigma(a_3)}q_{_{21}}^{-\tau(a_1-1)}Dy=q_{_{21}}^{-\tau(a_1)}yD$.
So {\small $$\begin{array}{ccl}
A&=&((-1)^{a_1-1+a_2}q_{_{12}}^{\sigma(a_3)\tau(a_1-1)}q_{_{21}}^{-\sigma(a_3)\tau(a_2)}[x_1^{\sigma(a_1)},x_2^{\sigma(a_3)}]_c\\
&&\hspace{5mm}+(-1)^{a_1}q_{_{21}}^{-\tau(a_1)}yD)\Phi(a_1-1,a_2,a_3-1)\\
&=&((-1)^{a_1-1+a_2}q_{_{12}}^{\sigma(a_3)\tau(a_1-1)}q_{_{21}}^{-\sigma(a_3)\tau(a_2)}[x_1^{\sigma(a_1)},x_2^{\sigma(a_3)}]_c\\
&&\hspace{5mm}+(-1)^{a_1}q_{_{21}}^{-\sigma(a_3)}q_{_{21}}^{-\tau(a_1-1)}Dy)\Phi(a_1-1,a_2,a_3-1)\\
&=&0.\end{array}$$}

The equations $B= 0$ and $C= 0$ can be verified directly. For
example,
{\small $$\begin{array}{ccl}B&=&(\dt_2\dt_3+\dt_3\dt_2)(\Phi(a_1,a_2, a_3))\\
&=&((-1)^{a_1}q_{_{21}}^{-\sigma(a_2)\tau(a_1)}y^{\sigma(a_2)}(-1)^{a_1+a_2-1}q_{_{12}}^{\sigma(a_3)\tau(a_1)}q_{_{21}}^{-\sigma(a_3)\tau(a_2-1)}x_2^{\sigma(a_3)}\\
&&\hspace{5mm}+(-1)^{a_1+a_2}q_{_{12}}^{\sigma(a_3)\tau(a_1)}q_{_{21}}^{-\sigma(a_3)\tau(a_2)}x_2^{\sigma(a_3)}(-1)^{a_1}q_{_{21}}^{-\sigma(a_2)\tau(a_1)}y^{\sigma(a_2)})\\&&\hspace{5mm}\Phi(a_1,a_2-1,a_3-1)\\
&=&0,\end{array}$$} since $\tau(a_2-1)+\sigma(a_2)=\tau(a_2)$.

If $a_2$ is even, then {\small $$\begin{array}{ccl}
\partial^2(\Phi(a_1,a_2, a_3))&=&((\dt_1\dt_3+\dt_3\dt_1+\dt_2\tilde{\dt}_2)+(\dt_1\dt_2+\dt_2\dt_1)+(\dt_3\dt_2+\dt_2\dt_3)\\
&&\hspace{5mm}+(\tilde{\dt}_2\dt_1+\dt_1\tilde{\dt}_2)+(\tilde{\dt}_2\dt_3+\dt_3\tilde{\dt}_2))\Phi(a_1,a_2,
a_3).\end{array}$$} The equation
$(\dt_1\dt_3+\dt_3\dt_1+\dt_2\tilde{\dt}_2)\Phi(a_1,a_2, a_3)=0$
follows directly from the definition of $\tilde{\dt}_2$. As in the
case in which $a_2$ is odd, $$(\dt_2\dt_3+\dt_3\dt_2)\Phi(a_1,a_2,
a_3)=0 \t{ and }(\dt_1\dt_2+\dt_2\dt_1)\Phi(a_1,a_2, a_3)=0$$ can be
also verified via a straightforward computation. Now, we show that
$(\tilde{\dt}_2\dt_1+\dt_1\tilde{\dt}_2)\Phi(a_1,a_2, a_3)=0$ case
by case, using Lemma \ref{d2}.

Case (i) $a_1$ and $a_3$ are both odd, {\small $$\begin{array}{ccl}
&&(\tilde{\dt}_2\dt_1+\dt_1\tilde{\dt}_2)\Phi(a_1,a_2,
a_3)\\&=&(x_1(q_{_{21}}^{-\frac{a_2}{2}}x_1^{N-2})-q_{_{21}}^{-\frac{a_2}{2}}x_1^{N-1})\Phi(a_1-2,a_2+1,a_3)\\
&=&0.\end{array}$$}

Case (ii) $a_1$ is odd and $a_3$ is even,
{\small$$\begin{array}{cl}&(\tilde{\dt}_2\dt_1+\dt_1\tilde{\dt}_2)\Phi(a_1,a_2, a_3)\\
=&(x_1(-q_{_{12}}^{(N-1)(\frac{a_1-3}{2}N+1)}q_{_{21}}^{-(N-1)\frac{a_2}{2}N}q_{_{21}}^{\frac{a_1-3}{2}N+1})(k_1x_1^{N-2}x_2^{N-2}+\cdots\\
&\hspace{5mm}+k_{N-2}y^{N-3}x_1x_2+k_{N-1}y^{N-2})\\
&\hspace{5mm}+q_{_{12}}^{(N-1)\frac{a_1-1}{2}N}q_{_{21}}^{-(N-1)\frac{a_2}{2}N}\bar{q}q_{_{21}}^{-(N-2)}q_{_{21}}^{\frac{a_1-1}{2}N}x_2^{N-2}x_1^{N-1})\Phi(a_1-2,a_2+1,a_3)\\
=&q_{_{12}}^{(N-1)\frac{a_1-1}{2}N}q_{_{21}}^{-(N-1)\frac{a_2}{2}N}\bar{q}q_{_{21}}^{-(N-2)}q_{_{21}}^{\frac{a_1-1}{2}N}\\
&\hspace{5mm}(-q_{_{12}}^{-N^2+2N}x_1(k_1x_1^{N-2}x_2^{N-2}+\cdots+k_{N-2}y^{N-3}x_1x_2+k_{N-1}y^{N-2})\\
&\hspace{5mm}+x_2^{N-2}x_1^{N-1})\Phi(a_1-2,a_2+1,a_3)\\
=&0,\end{array}$$} since
$q_{_{12}}^{-N^2+2N}x_1[x_1^{N-1},x_2^{N-1}]_c=x_2^{N-2}x_1^{N-1}y$.
% and the set $\{x_1^{N-i}x_2^{N-i-1}y^i\}$, $1\se i\se N-1$, is linearly independent.

Case (iii) $a_1$ is even and $a_3$ is odd,
{\small$$\begin{array}{cl}&(\tilde{\dt}_2\dt_1+\dt_1\tilde{\dt}_2)\Phi(a_1,a_2,
a_3)\\=&-q_{_{21}}^{-\frac{a_2}{2}N}x_1^{N-1}+q_{_{21}}^{-\frac{a_2}{2}N}x_1^{N-1}\Phi(a_1-2,a_2+1,a_3)\\=&0.\end{array}$$}

Case (iv) $a_1$ and $a_3$ are both even,
{\small$$\begin{array}{cl}&(\tilde{\dt}_2d_1+d_1\tilde{d}_2)\Phi(a_1,a_2, a_3)\\
=&x_1^{N-1}(q_{_{12}}^{(N-1)\frac{a_1-2}{2}N}q_{_{21}}^{-(N-1)\frac{a_2}{2}N}\bar{q}q_{_{21}}^{-(N-2)}q_{_{21}}^{\frac{a_1-2}{2}N}x_2^{N-2})\\
&\hspace{5mm}+(-q_{_{12}}^{(N-1)(\frac{a_1-2}{2}N+1)}q_{_{21}}^{-(N-1)\frac{a_2}{2}N}q_{_{21}}^{\frac{a_1-2}{2}N+1})(k_1x_1^{N-2}x_2^{N-2}+\cdots\\
&\hspace{5mm}+k_{N-2}y^{N-3}x_1x_2+k_{N-1}y^{N-2})x_1\Phi(a_1-2,a_2+1,a_3)\\
=&(q_{_{12}}^{(N-1)(\frac{a_1-2}{2}N+1)}q_{_{21}}^{-(N-1)\frac{a_2}{2}N}q_{_{21}}^{\frac{a_1-2}{2}N+1})(x_1^{N-1}x_2^{N-2}\\
&\hspace{5mm}-(k_1x_1^{N-2}x_2^{N-2}+\cdots+k_{N-2}y^{N-3}x_1x_2+k_{N-1}y^{N-2})x_1)\Phi(a_1-2,a_2+1,a_3)\\
=&0,\end{array}$$} since
$[x_1^{N-1},x_2^{N-1}]_cx_1=yx_1^{N-1}x_2^{N-2}$.

Similarly, we can prove that
$(\tilde{\dt}_2\dt_3+\dt_3\tilde{\dt}_2)\Phi(a_1,a_2,a_3)=0$.

In conclusion, we have $\partial^2=0$.

\subsection{}\label{exta2. append relation}

In this subsection, we give the necessary commutative diagrams to
check the relations in Theorems \ref{1} and \ref{2}.

Set {\small
$$X_1=q_{_{12}}^{-(N-1)(N-3)}x_1^{N-3}x_2^{N-3}+k_1yx_1^{N-4}x_2^{N-4}+\cdots+k_{N-3}y^{N-3},$$}
where $k_i\in \kk$, $1\se i\se N-3$,  such that
$x_2^{N-1}x_1^{N-3}=X_1x_2^2$, and {\small
$$X_2=q_{_{12}}^{(N-3)(N-1)}x_2^{N-3}x_1^{N-3}+l_1yx_2^{N-4}x_1^{N-4}+l_2y^2x_2^{N-5}x_1^{N-5}+\cdots+l_{N-3}y^{N-3},$$}
where $l_i\in \kk$, $1\se i\se N-3$,  such that
$x_1^{N-1}x_2^{N-3}=X_2x_1^2$.

Let $f_1^i$, $f_2^i$ , $f_3^i$ and $g_1^j$, $g_2^j$ , $g_3^j$, $1\se
i\se 5$ and $j=1,2$ be the morphisms described  by the following
matrices:

$f_1^i$ is the $5\times 1$ matrix with 1 in the $i$-th position and
0 elsewhere,

%$$f_1^i=\left(\begin{array} {c}0\\1\\0\\0\\0
%\end{array}\right),$$

$$f_2^1={\tiny\left(\begin{array} {cc}1&0\\
0&q_{_{12}}^N\\
0&0\\
0&0\\
0&0\\
0&0\\
0&0
\end{array}\right)},f_2^2={\tiny\left(\begin{array} {cc}0&0\\
x_1^{N-3}&0\\
0&0\\
0&1\\
q_{_{12}}q_{_{21}}^{1-N}y^{N-2}x_2&-q_{_{21}}^{1-N}y^{N-2}x_1\\
0&0\\
0&0
\end{array}\right)},$$
$$f_2^3={\tiny\left(\begin{array} {cc}0&0\\
0&0\\
0&1\\
0&0\\
1&0\\
0&0\\
0&0
\end{array}\right)},f_2^4={\tiny\left(\begin{array} {cc}0&0\\
0&0\\
-q_{_{12}}^2q_{_{21}}^{N-1}y^{N-2}x_2&q_{_{12}}q_{_{21}}^{N-1}y^{N-2}x_1\\
1&0\\
0&0\\
0&q_{_{12}}^Nx_2^{N-3}\\
0&0
\end{array}\right)},f_2^5={\tiny\left(\begin{array} {cc}0&0\\
0&0\\
0&0\\
0&0\\
0&0\\
1&0\\
0&1
\end{array}\right)},$$

$$f_3^1={\tiny\left(\begin{array} {ccccc}1&0&0&0&0\\
0&q_{_{12}}^N&0&0&0\\
0&0&q_{_{12}}^N&0&0\\
0&0&0&q_{_{12}}^N&0\\
0&0&0&0&q_{_{12}}^N\\
0&0&0&0&0\\
0&0&0&0&0\\
0&0&0&0&0\\
0&0&0&0&0\\
0&0&0&0&0\\
0&0&0&0&0\\
0&0&0&0&0\\
\end{array}\right)},f_3^3={\tiny\left(\begin{array} {ccccc}0&0&0&0&0\\
0&0&0&0&0\\
q_{_{12}}^{-N^2+N}&0&0&0&0\\
0&0&0&0&0\\
0&0&0&0&0\\
0&1&0&0&0\\
0&0&1&0&0\\
0&0&0&1&0\\
0&0&0&0&0\\
0&0&0&0&q_{_{12}}^{N^2}\\
0&0&0&0&0\\
0&0&0&0&0\\
\end{array}\right)},$$
$$f_3^2={\tiny\left(\begin{array} {ccccc}0&0&0&0&0\\
1&0&0&0&0\\
q_{_{21}}^{-1}y^{N-2}x_2&0&0&0&0\\
0&q_{_{12}}^{-N}x_1^{N-3}&0&0&0\\
0&0&0&X_1&0\\
0&q_{_{21}}^{1-N}q_{_{12}}^2y^{N-2}x_2&q_{_{12}}^{N}&0&0\\
0&0&0&0&0\\
0&0&0&q_{_{12}}q_{_{21}}^{N-3}y^{N-2}x_2&0\\
0&0&0&0&q_{_{12}}^{2N}\\
0&0&0&0&0\\
0&0&0&0&0\\
0&0&0&0&0\\
\end{array}\right)},$$

$$f_3^4={\tiny\left(\begin{array} {ccccc}0&0&0&0&0\\
0&0&0&0&0\\
0&0&0&0&0\\
q_{_{12}}^{-N}&0&0&0&0\\
0&q_{_{12}}^{-N^2+2N}X_2&0&0&0\\
0&q_{_{12}}q_{_{21}}^{-N+3}y^{N-2}x_1&0&0&0\\
0&0&0&0&0\\
0&0&q_{_{21}}^N&q_{_{21}}^{N-1}y^{N-2}x_1&0\\
0&0&0&q_{_{12}}^{2N}x_2^{N-3}&0\\
0&0&0&0&q_{_{21}}^{-N+1}q_{_{12}}^2y^{N-2}x_1\\
0&0&0&0&q_{_{12}}^N\\
0&0&0&0&0
\end{array}\right)},$$

$$f_3^5={\tiny\left(\begin{array} {ccccc}0&0&0&0&0\\
0&0&0&0&0\\
0&0&0&0&0\\
0&0&0&0&0\\
q_{_{12}}^{-N^2+N}&0&0&0&0\\
0&0&0&0&0\\
0&0&0&0&0\\
0&0&0&0&0\\
0&1&0&0&0\\
0&0&1&0&0\\
0&0&0&1&0\\
0&0&0&0&1\\
\end{array}\right)},$$

$$g_1^1={\tiny\left(\begin{array}{c}1\\
0\end{array}\right)},g_1^2={\tiny\left(\begin{array}{c}0\\
1\end{array}\right)},$$

$$g_2^1={\tiny\left(\begin{array}{cc}x_1^{N-2}&0\\
\bar{q}q_{_{12}}^2x_2&(-q_{_{12}}-\bar{q}q_{_{12}})x_1\\
0&-q_{_{12}}y^{N-1}\\
0&x_2\\
0&0
\end{array}\right)},g_3^1={\tiny\left(\begin{array}{ccccc}
1&0&0&0&0\\
0&(-q_{_{12}}-\bar{q}q_{_{12}})x_1^{N-3}&0&0&0\\
0&0&0&0&0\\
0&0&0&\bar{q}q_{_{12}}^2&0\\
0&0&q_{_{12}}^{N}&0&0\\
0&0&0&0&q_{_{12}}^N\\
0&0&0&0&0\\
\end{array}\right)},$$

$$ g_2^2={\tiny\left(\begin{array}{cc}
0&0\\
x_1&0\\
y^{N-1}&0\\
(-q_{_{21}}-\bar{q}q_{_{21}})x_2&\bar{q}q_{_{21}}^2x_1\\
0&x_2^{N-2}\\
\end{array}\right)}, g_3^2={\tiny\left(\begin{array}{ccccc}
0&0&0&0&0\\
1&0&0&0&0\\
0&0&q_{_{21}}^{N}&0&0\\
0&\bar{q}q_{_{21}}^2&0&0&0\\
0&0&0&0&0\\
0&0&0&q_{_{12}}^N(-q_{_{21}}-\bar{q}q_{_{21}})x_2^{N-3}&0\\
0&0&0&0&1\\
\end{array}\right)}.$$

Then we have the following commutative diagrams
{\begin{equation}\label{f}\xymatrix{
   R^{12}\ar[r]\ar[d]_{f_3^i}&R^{7}\ar[r]\ar[d]_{f_2^i}&R^{5}\ar[d]_{f_1^i}\ar[r]\ar[dr]&R^2\ar[r]&R\\
   R^5\ar[r]&R^2\ar[r]&R\ar[r]&\kk&},\end{equation}}
{\begin{equation}\label{g}\xymatrix{
   R^{7}\ar[r]\ar[d]_{g_3^i}&R^{5}\ar[d]_{g_2^i}\ar[r]&R^2\ar[d]_{g_1^i}\ar[r]\ar[dr]&R\ar[r]&\kk\\
   R^5\ar[r]&R^2\ar[r]&R\ar[r]&\kk&}.\end{equation}}

It is also routine to check the commutativity of the diagrams
(\ref{f}) and (\ref{g}). But we need to mention that the  following
equations hold
$$X_1(\bar{q}q_{_{21}}^2x_1x_2-(q_{_{21}}+\bar{q}q_{_{21}})x_2x_1)=-q_{_{12}}^{-N^2+2N}\overline{D},$$
$$X_2(\bar{q}q_{_{12}}^2x_2x_1-(q_{_{12}}+\bar{q}q_{_{12}})x_1x_2)=-\overline{D},$$
which follow from Lemma \ref{bas} and the following two equations
$$\begin{array}{ccl}q_{_{12}}^{-N^2+2N}\overline{D}x_2&=&x_2^{N-1}x_1^{N-2}\\&=&X_1x_2^2x_1\\
&=&X_1(-\bar{q}q_{_{21}}^2x_1x_2+(q_{_{21}}+\bar{q}q_{_{21}})x_2x_1)x_2,\end{array}$$
$$\begin{array}{ccl}\overline{D}x_1&=&x_1^{N-1}x_2^{N-2}\\&=&X_2x_1^2x_2\\
&=&X_2(-\bar{q}q_{_{12}}^2x_2x_1+(q_{_{12}}+\bar{q}q_{_{12}})x_1x_2)x_1.\end{array}$$

\bibliography{}

\end{document}